\documentclass[12pt,notitlepage]{amsart} 
\usepackage{amssymb,array,delarray,epsfig,latexsym,rotating} 

\newlength{\sh}
\newlength{\fw}
\newlength{\jmr}
\newlength{\jfc}
\newlength{\bernd}
\newlength{\ioan}
\newlength{\agk}
\newlength{\gkz}
\newtheorem{cor}{Corollary}
\newtheorem{lemma}{Lemma}

\newtheorem{dfn}{Definition}

\newtheorem{van}{Vanishing Theorem for Resultants} 

\newtheorem{main}{Main Theorem} 
\newtheorem{thm}[main]{Theorem}
\newtheorem{rem}{Remark}	
\newtheorem{ex}{Example}

\newcommand{\eps}{\varepsilon}
\newcommand{\cA}{\mathcal{A}}
\newcommand{\cD}{\mathcal{D}}
\newcommand{\cE}{\mathcal{E}}

\newcommand{\oE}{{\bar{E}}}
\newcommand{\oF}{{\bar{F}}}
\newcommand{\oP}{{\bar{P}}}

\newcommand{\cH}{\mathcal{H}}

\newcommand{\K}{{\mathbb{K}}}

\newcommand{\cO}{\mathcal{O}}

\newcommand{\cP}{\mathcal{P}}

\newcommand{\supp}{\mathrm{Supp}}
\newcommand{\conv}{\mathrm{Conv}}
\newcommand{\aff}{\mathrm{AffRes}}

\newcommand{\thth}{{\underline{\mathrm{th}}}}

\newcommand{\nd}{{\underline{\mathrm{nd}}}}
\newcommand{\Pro}{\mathbb{P}}

\newcommand{\Pa}{\Pro^{\#A-1}_\K}
\newcommand{\Pn}{\Pro^n_\K}
\newcommand{\Q}{\mathbb{Q}}
\newcommand{\R}{\mathbb{R}}
\newcommand{\C}{\mathbb{C}}
\newcommand{\N}{\mathbb{N}}
\newcommand{\Z}{\mathbb{Z}}
 
\newcommand{\fii}{\varphi}

\newcommand{\ch}{\mathrm{char}}

\newcommand{\area}{\mathrm{Area}}

\newcommand{\divisor}{\mathrm{Div}}

\newcommand{\choo}{\mathrm{Chow}}
\newcommand{\pert}{\mathrm{Pert}}

\newcommand{\fan}{\mathrm{Fan}}

\newcommand{\res}{\mathrm{Res}}

\newcommand{\Sn}{\mathcal{S}^{n-1}}
\newcommand{\Non}{(\N\cup\{0\})^n}
\newcommand{\Zn}{\Z^n}
 
\newcommand{\Qn}{\Q^n} 
\newcommand{\Rn}{\R^n}
\newcommand{\Qns}{\Q^n\!\setminus\!\{\bO\}}

\newcommand{\Ks}{\K^*}
\newcommand{\Kn}{\K^n}
 
\newcommand{\Rsn}{{(\R^*)}^n}
\newcommand{\Csn}{{(\C^*)}^n}
\newcommand{\Ksn}{{(\Ks)}^n}
\renewcommand{\qed}{$\blacksquare$}
\newcommand{\cM}{\mathcal{M}}
\newcommand{\sM}{\mathcal{S\!M}}
\newcommand{\aM}{\cM^{\mathrm{ave}}}
\newcommand{\cR}{\mathcal{R}}

\newcommand{\cZ}{\mathcal{Z}}
\newcommand{\wP}{{\oP}} 
\newcommand{\bZ}{\cZ^{\times}} 
\newcommand{\hZ}{\bZ_{A,+}} 
\newcommand{\dZ}{\bZ_A}

\newcommand{\cC}{\mathcal{C}}

\newcommand{\cT}{\mathcal{T}}
\newcommand{\bO}{\mathbf{O}}
\newcommand{\V}{\mathrm{Vol}}

\renewenvironment{Bmatrix}{\begin{array}{*{20}{c}}}{\end{array}}

\begin{document}

\title[Solving Degenerate Sparse Polynomial Systems Faster]{Solving 
Degenerate Sparse Polynomial Systems Faster} 

\author{J. Maurice Rojas}\thanks{
Submitted for publication.  
Formerly titled ``Twisted Chow Forms and Toric Perturbations for 
Degenerate Polynomial Systems.'' This research was partially funded by an 
N.S.F. Mathematical Sciences Postdoctoral Fellowship.}

\address{Department of Mathematics, City University of Hong Kong, 
83 Tat Chee Avenue, Kowloon, HONG KONG} 

\email{mamrojas@math.cityu.edu.hk\\ {\it Web-Page:} 
http://www.cityu.edu.hk/ma/staff/rojas}

\dedicatory{This paper is dedicated to my son, Victor Lorenzo.}

\date{September 14, 1998} 

\begin{abstract}
Consider a system $F$ of $n$ polynomial equations in 
$n$ unknowns, over an algebraically closed field of arbitrary 
characteristic. 
We present a fast method to find a point in every irreducible 
component of the zero set $\cZ$ of $F$. Our techniques allow 
us to sharpen and lower prior complexity bounds for 
this problem by fully taking into account the monomial term structure. 
As a corollary of our development we also obtain new explicit formulae 
for the exact number of isolated roots of $F$ and the intersection 
multiplicity of the positive-dimensional part of $\cZ$. Finally, 
we present a combinatorial construction of non-degenerate polynomial  systems, with specified monomial term structure and maximally many 
isolated roots, which may be of independent interest. 
\end{abstract} 

\mbox{}\\
\vspace{-.5in} 
\maketitle 

\section{Introduction}
\label{sec:intro}
The rebirth of resultants, especially through the {\bf  
toric}\footnote{Other commonly used prefixes for this 
modern generalization of the classical resultant \cite{vdv} include: 
sparse, mixed, sparse mixed, $\cA$-, $(\cA_1,\ldots,\cA_k)$-, and Newton.} 
resultant \cite{gkz94}, has begun to provide a much needed alternative to 
Gr\"obner basis methods for solving polynomial systems. Continuing 
this philosophy, we will use toric geometry to derive 
significant speed-ups and extensions of resultant-based methods for 
solving polynomial systems with infinitely many roots. 

The importance of dealing with degenerate polynomial 
systems has been observed in earlier work on quantifier 
elimination over algebraically closed fields 
\cite{chigo,renegar,pspace}: 
Many reasonable algorithms for polynomial system solving fail 
catastrophically when presented with a system $F$ (of $n$ polynomials in $n$ 
unknowns) having a positive-dimensional zero set $\cZ$. Even worse, this kind 
of failure can also occur when $F$ has only finitely 
many roots, if $F$ has infinitely many roots ``at infinity.'' 
When such failures occur, it is of considerable benefit to the user to at 
least be given some sort of description of the zero-dimensional  
part of $\cZ$. 

We will present two new techniques for handling such degeneracies. 
The {\bf twisted Chow form} (cf.\ Main Theorem \ref{main:chow}) allows one to 
quickly coordinatize many (but not all) degenerate $\cZ$, simply by 
injecting some extra combinatorics into the classical 
$u$-resultant. Our second technique builds on the twisted Chow form and works 
for {\bf all} degenerate $\cZ$: The {\bf toric perturbation} (cf.\ Main Theorem 
\ref{main:pert}) refines and generalizes an earlier algebraic perturbation 
trick used by Chistov and Grigoriev \cite{chigo}, 
Renegar \cite{renegar}, and Canny \cite{gcp}. 

Our refinement takes sparsity into account and 
allows one to replace the polynomial degrees present in earlier 
complexity bounds by more intrinsic geometric parameters (cf.\ Main 
Theorems \ref{main:big} and \ref{main:pert}). We will see 
in sections \ref{sub:dense} and \ref{sec:complex} that our bounds are a 
definite improvement, sometimes 
even by a factor exponential in $n$. Our framework also allows us to 
work over any algebraically closed field (as opposed 
to some earlier restrictions to the complex numbers) and to 
isolate the zero-dimensional part of $\cZ$.  

We also derive four corollaries which may be of independent interest: 
\begin{enumerate}
\item{An explicit method to compute field extensions involving the 
roots of $F$ (Corollary \ref{cor:galois}).} 
\item{An explicit formula for the exact, as opposed to generic, 
number of isolated\footnote{By an {\bf isolated} root, we will simply 
mean a root not lying in a positive-dimensional component of 
$\cZ$.} roots of $F$ (Corollaries \ref{cor:count} and \ref{cor:double}).} 
\item{A combinatorial construction, within polynomial time 
for fixed $n$, of explicit $F$ with specified monomial 
term structure and no roots ``at infinity'' (Main Theorem 
\ref{main:fill}).}
\item{A lower bound (conjecturally an explicit formula) for the intersection multiplicity of the positive-dimensional part of $\cZ$ (Corollary \ref{cor:double}). }  
\end{enumerate} 

Our main results are stated precisely in section \ref{sec:sum}. 
We then give several simple examples of our main results in section 
\ref{sec:ex}. There we also give an intuitive discussion of roots 
``at infinity'' and show how our results include Canny's 
earlier {\bf generalized characteristic polynomial (GCP)} 
as a special case. Section \ref{sec:fill} then details our aforementioned 
combinatorial   construction of ``generic'' $F$ with specified monomial 
term structure. Our main results are then proved in section \ref{sec:chowtor}, 
and we discuss the computational complexity of our techniques in 
section \ref{sec:complex}.

\section{Summary of Main Results} 
\label{sec:sum} 
Before describing our results in detail, we will 
introduce some necessary notation:  
In what follows, we will let $\oF\!:=\!(f_1,\ldots,f_{n+1})$, where 
for all $i$, $f_i(x)\!=\!\sum_{a\in E_i}c_{i,a}x^a$, 
$E_i$ is a nonempty finite subset of $\Non$, and $x^a$ is
understood to be the monomial term $x^{a_1}_1\cdots x^{a_n}_n$. 
Given the $c_{i,a}$, we will be solving for $x\!:=\!(x_1,\ldots,x_n)$. 
So the $(n+1)$-tuple $\oE\!:=\!(E_1,\ldots,E_{n+1})$ thus 
controls which monomial terms are allowed to appear in our systems of 
equations. An accepted shorthand is to say that {\bf  
$\boldsymbol{\oF}$ is an $\boldsymbol{(n+1)\times n}$ polynomial system with 
support contained in $\boldsymbol{\oE}$}. (This generalizes in an obvious 
way to $k\times n$ systems.) 

Of course, our given polynomial systems 
will usually be $n\times n$, so we will let 
$F\!:=\!(f_1,\ldots,f_n)$ and $E\!:=\!(E_1,\ldots,E_n)$. 
We also let $\conv(B)$ denote the convex hull of (i.e., smallest convex 
set containing) a point set $B\!\subseteq\!\Rn$, and let 
$[k]\!:=\!\{1,\ldots,k\}$ for any positive integer $k$. 
An important geometric invariant for $n\times n$ systems of equations 
is $\boldsymbol{\cM(E)}$ --- the {\bf mixed volume} 
\cite{buza,schneider,gk94,isawres,ewald,mvcomplex} of the convex hulls 
of the $E_i$. For \mbox{$(n+1)\times n$} systems, we also have the following 
two important complexity-theoretic parameters: 
$\boldsymbol{R(\oE)}\!:=\!\sum^{n+1}_{i=1}\cM(E_1,\ldots,E_{i-1},E_{i+1},
\ldots,E_{n+1})$ and $\boldsymbol{S(\oE)}\!=\!\cO(\sqrt{n}e^n\aM_{\oE})$, 
where $\aM_{\oE}$ is the average value of $\cM(\cE)$ as $\cE$ ranges over 
all $n$-tuples $(\cE_1,\ldots,\cE_n)$ with $\cE_j\!\in\{E_1,\ldots,E_{n+1}\}$ 
for all $j\!\in\![n]$. The true definition of $S(\oE)$ depends on the 
efficiency of a particular class of algorithms described later in sections 
\ref{sub:pert}, \ref{sub:solve}, and \ref{sec:complex}.  

We will usually take all polynomial coefficients to be 
constants in a fixed algebraically closed field $\K$ or polynomials in 
$\K[s]$ for some new parameter $s$. Also, we let 
$\Ks\!:=\!\K\!\setminus\!\{0\}$ and 
$\Delta\!:=\!\conv(\{\bO,\hat{e}_1,\ldots,\hat{e}_n\})$, where 
$\bO\!\in\!\Rn$ denotes the origin and $\hat{e}_i\!\in\!\Rn$ is the 
$i^\thth$ standard basis vector. Finally, using $\#$ for set 
cardinality, let $\fii_A : \Ksn 
\longrightarrow \Pro^{\#A-1}_\K$ be the rational map defined by 
$x \mapsto [x^a \; | \;  a\!\in\!A]$. On occasion, we will extend 
the domain of $\fii_A$ to a suitable toric variety (cf.\ 
section \ref{sec:chowtor}). 

\subsection{Finding Points in All Components in Intrinsic Polynomial 
Time}\mbox{}\\ 
Our first main result allows us 
to efficiently use exact arithmetic to find a point in 
every irreducible component of $\cZ$. In what follows, $\cO^*(T)$ 
means $\cO(T\log^r T)$ for some constant $r\!>\!0$. 
\begin{main}
\label{main:big} 
Let $F$ be an $n\times n$ polynomial system with support contained 
in $E$, assume $\cM(E)\!>\!0$, and set $E_{n+1}\!=\!A\!=\!\Delta\!\cap\!\Zn$. 
Also let $\fii_A(\cZ)$ be the zero set\footnote{Zero sets in projective space 
(and more general toric varieties) are 
defined in section \ref{sec:chowtor}. } of $F$ in $\Pn$. 
Then we can find univariate polynomials $h,h_1,\ldots,h_n$ 
with the following properties:
\begin{enumerate}
\addtocounter{enumi}{-1}
\item{The degrees of $h$ and $h_1,\ldots,h_n$ are all bounded 
above by $\cM(E)$.} 
\item{For any root $\theta$ of $h$, define 
$\gamma(\theta)\!:=\!(h_1(\theta),\ldots,h_n(\theta))$. Then 
$\gamma(\theta)\!\in\!\Ksn \Longrightarrow \gamma(\theta)$ is a 
root of $F$. } 
\item{There is at least one $\gamma(\theta)$ in every 
irreducible component of $\fii_A(\cZ)\cap\Ksn$. In particular, the set of 
points $\{\gamma(\theta)\}_{h(\theta)=0}$ contains all the isolated roots of 
$F$ in $\Ksn$. } 
\item{Let $K$ be $\Q(c_{i,a} \; | \;  i\!\in\![n],a\!\in\!E_i)$ 
or $(\Z/p\Z)(c_{i,a} \; | \;  i\!\in\![n],a\!\in\!E_i)$, according as 
$\ch \K$ is zero or a prime $p$. Then all the coefficients 
of $h,h_1,\ldots,h_n$ (and all intermediate calculations thereof) are in 
$K$, or a degree $\lceil 2\log_p((n+1)\cM(E))\rceil$
algebraic extension of $K$, according as $\ch \K$ is zero or $p$. } 
\end{enumerate}
Furthermore, we can find $h,h_1,\ldots,h_n$ deterministically within 
$\cO^*(n^4\cM(E)^3R(\oE)^2S(\oE)^{2.376})$ arithmetic steps 
and $\cO(nS(\oE)^2)$ space. Finally, at the expense of replacing 
$E$ by $\bO\cup E\!:=\!(\{\bO\}\cup E_1,\ldots,
\{\bO\}\cup E_n)$, we can ensure that $\{\gamma(\theta)\}_{h(\theta)=0}$ 
includes all the isolated roots of $F$ in $\Kn$ as well.
\end{main} 
\begin{rem}
\label{rem:speed} 
The above time bound can be reinterpreted as 
``near-heptic in the number of roots of a system closely related to 
$F$.'' Also, depending on the combinatorial data $E$ and 
the algebraic data $\ch \K$, the above complexity bounds  
can be lowered considerably, especially if randomization 
is allowed. These improvements are detailed further in section 
\ref{sec:complex}. In particular, Main Theorem \ref{main:big} 
already improves an earlier intrinsic complexity bound due to 
Giusti, Heintz, Morais, and Pardo \cite{joos}.\footnote{It should 
be noted that \cite{joos} also deals with the more general problem 
of complexity bounds for polynomial system solving in terms 
of arithmetic networks and straight-line programs. }   
\end{rem} 
\begin{rem}
\label{rem:matroid}
The assumption that $\cM(E)\!>\!0$ can actually be checked in 
polynomial time, via lemma \ref{lemma:ed} of 
section \ref{sec:fill}. Furthermore, if $\cM(E)\!=\!0$, then we can simply 
add $\leq\!n$ appropriately chosen points to $E$ (within the same 
asymptotic time bound) to make $\cM(E)$ positive.
\end{rem}

Our first main theorem thus removes a geometric/complexity-theoretic bottleneck from solving polynomial systems: For example, fast 
algorithms for finding approximations within $\eps\!>\!0$ of all the 
roots of 
$F$ in $\Csn$, within time $\cO^*(12^n\cM(E)^2\log\log\frac{1}{\eps})$ 
(neglecting some preprocessing), have recently been announced by 
Mourrain and Pan \cite{mp98}. However, their algorithm assumes that $\cZ$ is 
zero-dimensional and $\K\!=\!\C$. On the other hand, the algorithms of 
Canny from \cite{pspace,gcp} {\bf can} handle\footnote{That is, construct 
$h,h_1,\ldots,h_n$ as in Main Theorem 1.} positive-dimensional $\cZ$, but 
they assume $\K\!=\!\C$ and result in a Las Vegas complexity bound of 
$\cO^*(nd_\Pi$\begin{tiny}$\begin{pmatrix} d_\Sigma+1\\ n \end{pmatrix}$
           \end{tiny}$\!\!^{^4})$, 
where $d_\Pi$ and $d_\Sigma$ are respectively the product and sum of 
the total degrees of the $f_i$. We will see in sections \ref{sub:dense} and 
\ref{sec:complex} that 
our algorithm above is at least this fast, and is in fact frequently 
much faster. We also point out that when $\cZ$ is 
positive-dimensional, Gr\"obner basis techniques for solving 
$F$ suffer from a worst-case arithmetic complexity {\bf doubly} 
exponential in $n$ \cite{mm}. 

Main Theorem \ref{main:big} is also useful for certain rationality 
questions via the following corollary, proved in section \ref{sub:galois}. 
\begin{cor}
\label{cor:galois}
Following the notation of Main Theorem \ref{main:big}, 
suppose now that $\ch \K\!=\!0$ and $F$ has only finitely many roots in 
$\Ksn$. Let $g$ be the greatest common divisor of $h$ and $\prod^n_{i=1} 
h_i$. Then $K(\zeta_i \; | \; 
(\zeta_1,\ldots,\zeta_n)\!\in\!\Ksn \mathrm{ \ is \ a \ root \ of \ } 
F)$ is exactly the splitting field of $g$. 
\end{cor}   
 
To make Main Theorem \ref{main:big} more precise, we now outline its 
underlying toric geometric techniques.

\subsection{Main Geometric Results}\mbox{}\\  
First recall that there is a natural addition of 
point sets in $\Rn$ defined by $B+B'\!:=\!\{b+b' \; | \; b\!\in\!B, 
b'\!\in\!B' \}$.
In the notation of \cite{rio,bux}, we can associate to any 
$(n+1)$-tuple of point sets in $\Zn$, $\oE$, a {\bf toric resultant} 
$\res_\oE(\oF)$. This important operator is amply detailed in 
\cite{sparseelim,combiresult,gkz94,isawres,intro}, so 
let us state our first geometric construction. 
\begin{dfn} 
\label{dfn:main1}
Let $P\!:=\!\sum^n_{i=1}\conv(E_i)$ and $\oP\!:=\!P+\conv(E_{n+1})$. 
Also let $A\!\subset\!\Zn$ be any finite subset with at least two points and 
define $f_{n+1}(x)\!:=\!\sum_{a\in A} u_a x^a$ and 
$u\!:=\!(u_a \; | \; a\!\in\!A)$, 
where the $u_a$ are new parameters. We then call 
$\boldsymbol{\choo_A(u)}\!:=\!\res_{(E,A)}(F,f_{n+1})$ 
a {\bf twisted Chow form} of $F$. (Frequently, we will 
set $E_{n+1}\!=\!A$ and thus $\oP\!=\!P+\conv(A)$ as well.) 
\end{dfn}  
\noindent
Note that $\choo_A(u)$ will be a polynomial in the 
parameters $u_a$, encoding (in a manner to be described below) the  
roots of $F$. Twisted Chow forms are a generalization of the 
classical {\bf $\boldsymbol{u}$-resultant} \cite{vdv} since the latter 
simply corresponds to the case where we use the classical ``dense'' resultant 
and let $A\!=\!\Delta\cap\Zn$. For convenience, we will frequently 
respectively write $u_0$ and $u_i$ in place of $u_\bO$ and $u_{\hat{e}_i}$. 
\begin{ex}
Suppose we take $\ch \K\!\not\in\!\{2,3\}$, $n\!=\!2$, 
$E_1\!=\!E_2\!=\!2\Delta\cap\Z^2$, $A\!=\!\Delta
\cap\Z^2$, and $F\!=\!(1+2y-x^2+y^2,1+2x+x^2-4y^2)$. 
Then $\choo_A$ is simply the $u$-resultant, and this polynomial 
in $u_0,u_1,u_2$ factors (modulo a nonzero constant multiple) 
as $(u_0+\frac{1}{3}u_1-\frac{2}{3}u_2)\times\\ 
(u_0+3u_1+2u_2)(u_0-u_1)^2$. 
It is also not hard to see that $F$ has exactly three roots: 
$(\frac{1}{3},-\frac{2}{3})$, $(3,2)$, and $(-1,0)$; 
the last occuring with multiplicity $2$. Better still, 
we can read this off directly from our $u$-resultant 
by computing $(\frac{\mathrm{coefficient \ of \ }u_1}
{\mathrm{coefficient \ of \ }u_0},
\frac{\mathrm{coefficient \ of \ }u_2}
{\mathrm{coefficient \ of \ }u_0})$ for each linear factor  
(with $u_0$ appearing) of the $u$-resultant. 
(See Main Theorem \ref{main:chow} below.)  
\end{ex}

Our next main theorem tells us exactly how and when 
we can use a twisted Chow form to compute monomials in the roots of 
$F$. Recall that to any $n$-dimensional rational polytope 
$Q\!\subset\!\Rn$ one can associate its corresponding {\bf toric variety} 
(over $\K$) $\cT(Q)$ \cite{kkms,dannie,ksz92,tfulton,gkz94,toricint}, and this 
$\cT(Q)$ always has\footnote{It is {\bf not} always the 
case that $\cT(Q)$ also has a naturally embedded copy of $\Kn$. 
However, with some extra work, one can modify $Q$ so that this is true.} a 
naturally embedded copy of $\Ksn$. To state our results fully, we will 
require some toric variety terminology, but the underlying idea is 
simple: By working in compactifications more general than the 
projective spaces $\{\Pn\}^\infty_{n=1}$, we can make better use of the 
monomial term structure of our polynomial systems. 
\begin{main}
\label{main:chow} 
Following the notation of definition \ref{dfn:main1}, 
set $E_{n+1}\!=\!A$
and let $\cZ$ denote the zero set of $F$ in $\cT(\oP)$. Then 
$\choo_A(u)$ is a homogeneous polynomial, either identically 
zero or of degree $\cM(E)$, with the 
following properties:
\begin{enumerate}
\item{The polynomial $\choo_A$ is indentically zero  
$\Longleftrightarrow \fii_A(\cZ)$ is positive-dimensional. } 
\item{If $\zeta\!\in\!\cT(\oP)$ is a root of $F$ then 
$\choo_A$ is divisible by $\sum_{a\in A} \gamma_a u_a$, where 
\mbox{$[\gamma_a\; |\; a\!\in\!A]\!=\!\fii_A(\zeta)$.}}  
\item{The polynomial $\choo_A(u)$ splits completely (over $\K$) 
into linear factors. In particular, if $\choo_A\!\not\equiv\!0$ and a 
nonzero linear form $\sum_{a\in A} \gamma_a u_a$ divides 
$\choo_A$, then $[\gamma_a \; | \; a\!\in\!A]\!=\!\fii_A(\zeta)$ for {\bf some} 
root $\zeta\!\in\!\cT(\oP)$ of $F$.}
\end{enumerate} 
\end{main} 
\noindent
The zero set of $F$ in a toric variety is formalized in section 
\ref{sec:chowtor}. 
Note in particular that assertions (2) and (3) 
tell us that calculating $\choo_A(u)$ allows us to 
reduce the computation of the projective coordinates  
$[\zeta^a \; | \; a\!\in\!A]$, for any root $\zeta\!\in\!\cT(\oP)$ of 
$F$, to a multivariate factorization problem. 
Of course, this reduction only works if $\choo_A(u)$ is 
not indentically zero, and assertion (1) tells us exactly 
when this happens. 

We also obtain the following almost immediate corollary. 
\begin{cor}
\label{cor:count} 
Following the notation of Main Theorem \ref{main:chow}, we may check if 
$\choo_A$ is identically zero (and thus whether $\dim \fii_A(\cZ)\!>\!0$) 
within $\cO^*(n^2\cM(E)R(\oE)S(\oE)^{2.376})$ 
arithmetic steps and $\cO(nS(\oE)^2)$ space.\footnote{Just 
as in Main Theorem 1, these complexity bounds can be significantly lowered 
under certain reasonable assumptions. Also, unless otherwise stated, arithmetic 
steps will always be counted over the finite extension of $K$ described in 
Main Theorem \ref{main:big}. } Furthermore, if $\choo_A(u)$ does not vanish 
identically, then we can compute the exact number of roots of $F$ in $\Ksn$, 
counting multiplicities, within $\cO^*(n^4\cM(E)^3R(\oE)S(\oE)^{2.376})$ arithmetic steps and $\cO(nS(\oE)^2)$ space.$^7$ 
\end{cor} 
\noindent 
Even better, by combining with corollary \ref{cor:parts} 
of section \ref{sec:chowtor}, we can also see how many roots lie at various 
parts of ``toric infinity.'' Corollary \ref{cor:count} thus generalizes 
Bernshtein's famous mixed volume bound \cite{bernie} to {\bf exact} root 
counting over an algebraically closed field. 

However, there is still another improvement to be made:  
It is actually possible for $F$ to have infinitely 
many roots in $\cT(\oP)$ but only finitely many roots in $\Ksn$. 
In such cases, {\bf sometimes} the right $A$ will permit an 
exact count of the roots of $F$ in $\Ksn$ via Corollary \ref{cor:count}.  
For example, it is easy to construct 
$F$, $A$, and $A'$ where $\choo_A$ vanishes identically but 
$\choo_{A'}$ does not (cf.\ section \ref{sub:stranger}). On the other hand, 
those $F$ with infinitely many roots in $\Ksn$ will never 
have a nontrivial twisted Chow form.  

Our next construction works for {\bf all} $F$ and $A$, and begins as 
follows:
\begin{dfn} 
\label{dfn:pert} 
Following the notation of Main Theorem \ref{main:chow}, assume further that 
\mbox{$\cM(E)\!>\!0$.} Let $F^*$ be any $n\times n$ system 
with constant coefficients and support contained in $E$, 
such that $F^*$ has only finitely many roots in $\cT(P)$. We then 
say that $\boldsymbol{\cH(u;s)}\!:=\!\res_{(E,A)}(F-sF^*,f_{n+1})$ 
(where $s$ is a new indeterminate) is a {\bf toric generalized 
characteristic polynomial for $\boldsymbol{(F,A)}$.} Furthermore, 
we define $\boldsymbol{\pert_{A,F^*}(u)}\!\in\!K[u_a \; | \; 
a\!\in\!A]$ to be the coefficient of the term of $\cH(u;s)$ of 
lowest degree in $s$. We call $\pert_{A,F^*}$ a 
{\bf toric perturbation of $\boldsymbol{(F,A)}$} and, when no 
confusion is possible, we will sometimes write $\boldsymbol{\pert_A}$ instead. 
\end{dfn} 
 
The polynomial $\pert_A$ is what we can use in place of 
$\choo_A$ when $\choo_A$ vanishes identically. We will 
describe this shortly, but first we digress momentarily 
to describe how to construct the necessary ``generic'' 
$F^*$ above: If we simply fix the support of $F^*$ to be $E$, and 
pick random numbers for the coefficients (using any probability distribution 
on $\K^{\#\mathrm{monomial \ terms}}$ yielding probability $1$ avoidance of 
algebraic hypersurfaces), lemma \ref{lemma:bernie} of section 
\ref{sec:chowtor} tells us that $F^*$ will satisfy the above hypothesis with 
probability $1$. Alternatively, a deterministic method for constructing 
suitable $F^*$ is the following. 
\begin{dfn} 
\label{dfn:fill}
\cite{convexapp,rojaswang}
Given $n$-tuples $D\!:=\!(D_1,\ldots,D_n)$ and $E\!:=\!(E_1,\ldots,E_n)$ 
of nonempty compact subsets of $\Rn$, we say that $D$ 
{\bf fills} $E$ (or $D$ is a fill of $E$) iff (0) $D_i\!\subseteq\!E_i$ for all 
$i\!\in\![n]$ and (1) $\cM(D)\!=\!\cM(E)$. 
We then call $D$ {\bf irreducible} iff 
the removal of any point of $D$ causes $\cM(D)$ to decrease. 
\end{dfn}
\begin{main}
\label{main:fill} 
Following the notation of definition \ref{dfn:fill}, 
suppose $E_i\!\subset\!\Zn$ for all $i$, $\cM(E)\!>\!0$, and $D$ is an 
irreducible fill of $E$. Then, for any choice of {\bf nonzero} 
$c_{i,a}\!\in\!\Ks$, the polynomial system $(\sum_{a\in D_1} 
c_{1,a}x^a,\ldots, \sum_{a\in D_n} c_{n,a}x^a)$ has exactly $\cM(E)$ 
roots, counting multiplicities, in $\Ksn$ and no roots 
in $\cT(P)\!\setminus\!\Ksn$. Furthermore, letting $m\!:=\!\sum^n_{i=1} \#E_i$, an irreducible fill of $E$ can be found within time 
$\cO(nm^{2n+c+1})$, for some absolute constant $c\!>\!0$. 
\end{main}
\noindent  
Some simple examples of fills appear in section \ref{sub:fill} and 
we present further background on filling in section \ref{sec:fill}.  
We emphasize that while it is much more practical to pick a generic 
$F^*$ via randomization, the cost of derandomizing via fills 
can sometimes be amortized when one solves many $F$ with 
similar monomial term structure. In particular, the selection of an $F^*$ 
need only be done {\bf once} for a given $n$-tuple $E$, regardless of 
the coefficients of $F$.  

Toric perturbations improve on twisted Chow forms as 
follows:
\begin{main}
\label{main:pert} 
Following the notation of definition 2, $\pert_A(u)$ is a {\bf nonzero} 
homogeneous polynomial of degree $\cM(E)$ with the following properties:
\begin{enumerate}
\item{$\choo_A\!\not\equiv\!0 \Longleftrightarrow \cH(s)$ has a  
nonzero constant term. Also, when the latter holds, $\choo_A\!=\!\pert_A$. }
\item{If $\zeta\!\in\!\cT(\oP)$ is an {\bf isolated} root of 
$F$ then $\pert_A$ is divisible by $\sum_{a\in A} \gamma_a u_a$, where 
\mbox{$[\gamma_a\; | \; a\!\in\!A]\!=\!$} $\fii_A(\zeta)$.}  
\item{The polynomial $\pert_A(u)$ splits completely (over $\K$) 
into linear factors. In particular, extending the correspondence  
of assertion (2), for every irreducible positive-dimensional 
component $W$ of $\cZ$, there is at least one factor of $\pert_A$ 
corresponding to a root $\zeta\!\in\!W$. } 
\end{enumerate} 
Furthermore, we may evaluate $\pert_A$ at any point in $\K^{\#A}$ 
within $\cO^*(nR(\oE)^2S(\oE)^{2.376})$ arithmetic steps over $\K$ and $\cO(nS(\oE)^2)$ space.\footnote{Just 
as in Corollary \ref{cor:count} and Main Theorem 1, these complexity 
bounds can also be significantly lowered under certain reasonable assumptions.}  
\end{main}
\noindent 
We emphasize that the main advantage of $\pert_A$ is that we can 
pick {\bf any} $A$ we prefer and still get a useful 
analogue of $\choo_A$. For instance, even if the $u$-resultant 
unluckily vanishes identically, we can always 
simply set $A\!=\!\Delta\cap\Zn$ and directly read off the coordinates 
of the isolated roots of $F$ from the factors of $\pert_A(u)$ 
(assuming one can do multivariate factoring over $\K$). 
Indeed, $\pert_{\Delta\cap\Zn}$ and assertion (3) are central 
to our explicit construction of points in every irreducible 
component\footnote{The analogue of assertion (3) had been conjectured for 
Canny's GCP. We have thus proved this conjecture and generalized 
it to the toric GCP.} of $\cZ$, not to mention the proof of Main Theorem 
\ref{main:big}. 

Better still, we can sometimes (conjecturally always) distinguish which roots of $F$ are isolated. 
\begin{cor} 
\label{cor:double} 
Following the notation above, let $\cZ_0$ and $\cZ_\infty$ 
respectively denote the zero-dimensional and positive-dimensional parts of $\cZ$. Then $\cZ_\infty\cap\Ksn\!=\!\emptyset 
\Longrightarrow$ we can count the number of points in 
$\cZ_0\cap\Ksn$, with or without multiplicity, within the same asymptotic complexity bounds as stated in Main Theorem \ref{main:big}. More generally, there is a randomized algorithm 
which computes upper bounds on the cycle class degrees $\deg \cZ_0$ and 
$\deg \cZ_0\cap\Ksn$, and a lower bound on $\deg \cZ_\infty$, within the same complexity bounds. Conjecturally, these bounds are 
all actually {\bf explicit formulae} with probability $1$. 
\end{cor} 
\noindent 
A simple example of this final main result (and Main Theorem 
\ref{main:pert}) also appears 
in section \ref{sub:pert}. So in summary, as the 
zero set of $F$ in $\cT(\oP)$ becomes more and more degenerate, we can successively use Corollaries \ref{cor:count} and \ref{cor:double}  
to count roots in $\Ksn$ with complete generality. 
We also point out that a special case of Corollary \ref{cor:double} was used in \cite{venice} in connection with a fast general algorithm for exact multivariate root counting in $\Rsn$. 

We can also construct the corresponding analogues of 
$h$ and the $h_i$ to describe $\cZ_0$ explicitly, but this becomes more 
technical (cf.\ section \ref{sub:double}). The 
same can be said for the analogous results in $\Kn$, and this is covered in greater depth in \cite{bux} and \cite{aff}. 
We thus obtain a first step toward an algorithmic foundation for excess 
intersections. (See \cite{ifulton} for a brief historical 
description of this problem.) 

We now illustrate our results and theory.

\section{Examples}  
\label{sec:ex} 
We begin with two small examples of filling.
We will then see applications of the toric GCP and twisted Chow form to some 
degenerate $2 \times 2$ and $3 \times 3$ polynomial systems. Finally, 
we will see a brief comparison of the toric GCP to the original GCP. In what  
follows, we will sometimes respectively write $x$, $y$, and $z$ in place of 
$x_1$, $x_2$, and $x_3$. 

\subsection{Filling Squares and Cubes}\mbox{}\\
\label{sub:fill}
For our first example, consider the pair of rectangles 
$\cP\!:=\!([0,a]\!\times\![0,b],[0,c]\!\times\![0,d])$ where
$a$, $b$, $c$, and $d$ are positive integers. Then it 
is easily verified (via theorem \ref{thm:me} of section 
\ref{sec:fill}) that the pair $D\!=\!(\{(0,0),(a,b)\},\{(0,d),(c,0)\})$
fills $\cP$. In this case, the mixed area of both pairs is easily 
checked to be $ad+bc$. Note also that $D$ is a pair of oppositely 
slanting diagonals of our initial pair of rectangles (modulo 
taking convex hulls). Finally, it is easily checked that $D$ is 
indeed irreducible, since the removal of any point of $D$ results 
in a mixed area of $0$. 

By Main Theorem \ref{main:fill}, we thus obtain that for {\bf any} 
$\alpha_1,\alpha_2,\beta_1,\beta_2\!\in\!\Ks$, the bivariate polynomial 
system $(\alpha_1+\alpha_2x^ay^b,\beta_1x^a+\beta_2y^b)$ 
will have exactly $ad+bc$ roots, counting multiplicities, in $(\Ks)^2$. 

For our second example, let $\cP$ instead be a triple of 
standard cubes (so that the vertex set of each cube is simply 
$\{0,1\}^3$). Then, using the criterion from theorem \ref{thm:me} 
once again, it is easily verified that the triple 
$D\!=\!(\{(1,0,0),(0,1,0),(0,0,1)\},\{(1,1,0),(1,0,1),(0,1,1)\},\\ 
\{(0,0,0),(1,1,1)\})$ 
fills $\cP$. (This is depicted in Figure 1 below.) Also, it is easily checked 
that the mixed volume of both triples is $6$. Finally, note that this 
$D$ is irreducible as well by theorem \ref{thm:me}. Alternatively, one 
can easily check this by brute force, using any one of the publically 
web-accessible software packages for mixed volume computation by Emiris, 
Gao, Huber, or Verschelde. 

\begin{figure}[h]
\begin{center}
\epsfig{height=1in,file=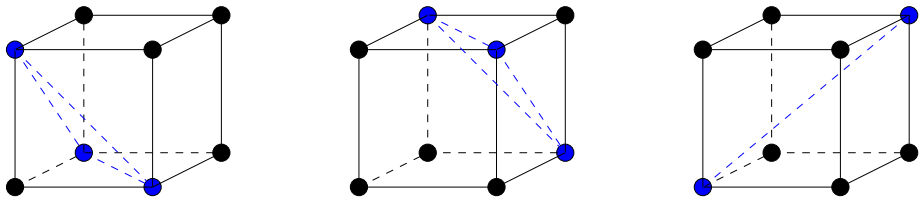}
\end{center} 
\caption{An irreducible fill of three $3$-cubes.} 
\end{figure}

By Main Theorem \ref{main:fill}, we thus obtain that for {\bf any} 
$\alpha_1,\alpha_2,\beta_1,\beta_2,\beta_3,\gamma_1,\gamma_2,
\gamma_3\!\in\!\Ks$, the trivariate polynomial system 
$(\alpha_1+\alpha_2xyz,\beta_1x+\beta_2y+\beta_3z,\gamma_1xy 
+\gamma_2xz+\gamma_yz)$ will have exactly $6$ roots, counting 
multiplicities, in $(\Ks)^3$. 

In summary, theorem \ref{thm:me} of section \ref{sec:fill} 
gives a necessary and sufficient criterion for $D$ to fill 
a given $n$-tuple $E$, and Main Theorem \ref{main:fill} tells us that we 
can construct {\bf some} irreducible fill for $E$ within time singly exponential in $n$.

\subsection{$\boldsymbol{\pert_A}$ Applied to a Degenerate 
$\boldsymbol{2\times 2}$ System}\mbox{}\\ 
\label{sub:pert} 
Consider the bivariate polynomial system 
\[ F\!=\!(1+2x-2x^2y-5xy+x^2+3x^3y,2+6x-6x^2y-11xy+4x^2+5x^3y)\] 
over any field of characteristic not equal to $2$, $3$, or $7$.  
Letting $E$ be the support of $F$, the reader can easily verify that 
$\cM(E)\!=\!4$, and that the only roots of $F$ are the points 
$\{(1,1),(\frac{1}{7},\frac{7}{4})\}$ and the line $\{-1\}\!\times\!\K$.\footnote{\mbox{For $n=2$, there is the simple 
formula $\cM(E)=\area(\conv(E_1+E_2))-\area(\conv(E_1))-\area(\conv(E_2))$.} 
Also, both polynomials are divisible by $x+1$. Furthermore, 
when $\ch \K\!=\!2$, the 
second isolated root becomes an isolated root lying on the $x$-axis.} So
it would appear that the $u$-resultant (and even $\choo_{\Delta\cap\Z^2}$) 
will vanish identically and not give us any useful information about any of 
these roots. Let us see how we can use $\pert_A$ (with $A\!=\!\Delta\cap\Z^2$) to recover everything we need to know about 
the roots of $F$. 

First, via combinatorial means \cite{sparseelim,isawres}, 
we construct a {\bf toric resultant matrix}, $M_\oE$. This matrix 
has the property that its determinant is a multiple of the toric resultant defining the toric GCP (the precursor to $\pert_A$). With the  
assistance of a {\tt Matlab} program, {\tt res2.m} (publically available from 
the author's web-page), we can obtain the following $17\times 17$ matrix: 
\footnotesize  
\[ M_\oE = \left[ \begin{Bmatrix} 

 u_1 & 0 & 0 & 0 & 0 & 0 & 0 & 0 & u_2 & 0 & 0 & 0 & 0 & 0 & 0 & 0 & u_0 \\ 

 u_0 & u_1 & 0 & 0 & 0 & 0 & 0 & 0 & 0 & u_2 & 0 & 0 & 0 & 0 & 0 & 0 & 0 \\ 
 
 0 & 0 & u_2 & 0 & 0 & 0 & 0 & 0 & u_1 & 0 & 0 & 0 & 0 & 0 & 0 & u_0 & 0 \\ 
 
 0 & 0 & 0 & u_2 & 0 & 0 & 0 & 0 & u_0 & u_1 & 0 & 0 & 0 & 0 & 0 & 0 & 0 \\ 
 
 0 & 0 & b_3 & b_4 & b_5 & 0 & 0 & 0 & b_2 & 0 & 0 & 0 & 0 & 0 & b_0 & b_1 & 0 \\ 
 
 0 & 0 & 0 & b_3 & b_4 & b_5 & 0 & 0 & b_1 & b_2 & 0 & 0 & 0 & 0 & 0 & b_0 & 0 \\ 
 
 0 & 0 & 0 & 0 & b_3 & b_4 & b_5 & 0 & b_0 & b_1 & 0 & 0 & b_2 & 0 & 0 & 0 & 0 \\ 
 
 0 & 0 & 0 & 0 & 0 & b_3 & b_4 & b_5 & 0 & b_0 & 0 & 0 & b_1 & b_2 & 0 & 0 & 0 \\ 
 
 0 & 0 & 0 & 0 & a_3 & a_4 & a_5 & 0 & a_0 & a_1 & 0 & 0 & a_2 & 0 & 0 & 0 & 0 \\ 
 
 0 & 0 & 0 & 0 & 0 & a_3 & a_4 & a_5 & 0 & a_0 & 0 & 0 & a_1 & a_2 & 0 & 0 & 0 \\ 
 
 0 & 0 & 0 & 0 & 0 & 0 & 0 & 0 & a_4 & a_5 & a_0 & a_1 & 0 & 0 & 0 & a_3 & a_2 \\ 

 a_2 & 0 & 0 & 0 & 0 & 0 & 0 & 0 & a_3 & a_4 & 0 & a_0 & a_5 & 0 & 0 & 0 & a_1 \\ 
 
 b_2 & 0 & 0 & 0 & 0 & 0 & 0 & 0 & b_3 & b_4 & 0 & b_0 & b_5 & 0 & 0 & 0 & b_1 \\ 
 
 b_1 & b_2 & 0 & 0 & 0 & 0 & 0 & 0 & 0 & b_3 & 0 & 0 & b_4 & b_5 & 0 & 0 & b_0 \\ 
 
 0 & 0 & a_3 & a_4 & a_5 & 0 & 0 & 0 & a_2 & 0 & 0 & 0 & 0 & 0 & a_0 & a_1 & 0 \\ 
 
 0 & 0 & 0 & a_3 & a_4 & a_5 & 0 & 0 & a_1 & a_2 & 0 & 0 & 0 & 0 & 0 & a_0 & 0 \\ 
 
 a_1 & a_2 & 0 & 0 & 0 & 0 & 0 & 0 & 0 & a_3 & 0 & 0 & a_4 & a_5 & 0 & 0 & a_0  
\end{Bmatrix} 
\right] \]
\normalsize  
where the $a_i$ (resp.\ $b_i$) are indeterminates correponding to the 
coefficients of $f_1$ (resp.\ $f_2$). Note in particular that $R(\oE)\!=\!4+4+4\!=\!12$. As for the other complexity parameter $S(\oE)$, its true definition is 
the size of any available toric resultant matrix. So $S(\oE)\!=\!17$ 
in the case at hand. 

Now note that by theorem \ref{thm:me}, 
$D\!:=\!(\{\bO,(3,1)\},\{(1,1),(2,0)\})$ is an irreducible fill of $E$. 
So by Main Theorem \ref{main:fill}, 
we can take $F^*\!=\!(1+x^3y,xy+x^2)$ and apply Main Theorem \ref{main:pert} to construct the toric GCP, $\cH(u;s)$. By setting 
$(a_0,\ldots,a_5)=(1-s,2,-2,-5,1,3-s)$, 
$(b_0,\ldots,b_5)\!=\!(2,6,-6,-11-s,4-s,5)$, and 
taking the determinant of $M_\oE$, we then obtain a nonzero constant 
multiple of $\cH(u_0,u_1,u_2;s)$.    

However, multivariate symbolic expansions are typically 
slow and memory-intensive. So to efficiently ``solve'' $F$ --- that is, 
to quickly find a point in every irreducible component of its zero set 
--- we will instead compute the univariate polynomials $h,h_1,h_2$ from 
Main Theorem 1 via interpolation. The polynomial $h$ is derived 
simply by specializing $\pert_A$ at some suitable value of $(u_1,u_2)$ 
and then interpolating through $1+\cM(E)$ values of $u_0$. The derivation 
of $h_1$ and $h_2$ is essentially the same but involves an additional intermediate step described in section \ref{sub:solve}. Since $\pert_A$ is in turn a coefficient of $\cH(u;s)$, there is also another level 
of interpolation through $1+S(\oE)-\cM(E)$ values of $s$.    

For example, setting $(u_1,u_2)\!=\!(\frac{1}{2},1)$ 
(and setting $u_0$ equal to a parameter $t$), we easily obtain via 
{\tt Maple} that  
\[h(t)  =  -153+120t+1540t^2+1600t^3+448t^4\] 
\[ h_1(t) = -\frac{11762}{7511}+\frac{19150}{22533}t
       +\frac{114736}{22533}t^2+\frac{7264}{3219}t^3\]
\[ h_2(t) = -\frac{5881}{7511}+\frac{32108}{22533}t
       +\frac{57368}{22533}t^2+\frac{3632}{3219}t^3.\]    
Since $h(t)$ factors as 
$(2t+1)(2t+3)(4t-1)(28t+51)$, we thus 
immediately obtain from Corollary \ref{cor:galois} (and the 
fact that $u_1$ and $u_2$ were chosen within $K$) that the 
zero-dimensional part of $\cZ\cap (\Ks)^2$ actually lies in $(K^*)^2$, 
where $K$ is the quotient field generated by the canonical image of $\Z$ 
in $\K$. Furthermore, by Main Theorem 1 (and the fact that 
$\bO\!\in\!E_1\cap E_2$), we can at last recover a set 
of points lying in $\cZ$ (including all 
the isolated roots of $F$ in $\K^2$) by substituting $\{-\frac{1}{2}, 
-\frac{2}{3},\frac{1}{4},-\frac{51}{28}\}$ into the pair 
$(h_1(t),h_2(t))$. 

As for the full expansion of the toric GCP, we can easily compute via {\tt Maple} that $\cH(u;s)$ is, up to a nonzero constant multiple,  
\scriptsize
\[ 
(u_2^4-u_0^4+u_1^4+6u_1^2u_2^2-4u_1u_2^3-4u_1^3u_2) \mathbf{s^8} 
\]
\[+(36u_1^2u_2^2-20u_2u_0^3-20u_2^3u_0- 4u_1u_0^3-19u_0^4-24u_2^4 
+6u_0^2u_1u_2 \]
\[ +36u_1u_2^3+36u_1^4-12u_0u_1
^2u_2-9u_1^2u_0^2+3u_2^2u_0^2+36u_0u_1u_2^2-4u_0u_1^3-84u_1^3u_2)\mathbf{s^7}\]
\[ +(220u_2^4-170u_2
u_0^3-394u_1^3u_2-98u_1u_0^3-98u_0^2u_1u_2-20u_0^4+370u_2^3u_0\] 
\[ +14u_0u_1
u_2^2-110u_0u_1^3-226u_1^2u_0^2-354u_1^2u_2^2+454u_1^4-274u_0u_1^2u_2+
74u_1u_2^3)\mathbf{s^6}\]
\[ +(1008u_2u_0^3-1612u_0^2u_1u_2+903u_0^4-624u_1u_0^3-2632u_2^3u_0-2104u_0
u_1^2u_2-970u_2^4\]
\[ -1010u_1u_2^3 +418u_1^3u_2 -2104u_0u_1u_2^2-642u_1^2u_2^2-1547
u_1^2u_0^2-936u_0u_1^3-1557u_2^2u_0^2+2204u_1^4)\mathbf{s^5}\]
\[ +(538u_0^2u_1u_2+1271u_0^4+
12253u_2^2u_0^2+6972u_2u_0^3+1929u_1^4-3075u_1^2u_2^2+654u_0u_1u_2^2\] 
\[+50u_1u_0^3 +2156u_2^4-960u_1^2u_0^2-2290u_0u_1^3+132u_1u_2^3-
5344u_0u_1^2u_2-1142u_1^3u_2+8708u_2^3u_0)\mathbf{s^4} \] 
\[ +(4384u_1u_0^3-24988u_2^2u_0^2-1582u_1^3u_2-6756u_0^4+10884u_0
u_1u_2^2+3802u_1u_2^3+15438u_0u_1^2u_2\] 
\[ 
+1024u_0u_1^3+8324u_1^2u_0^2-12826u_2^3
u_0+11270u_0^2u_1u_2-6976u_1^4+7164u_1^2u_2^2-21326u_2u_0^3-2408u_2^4)
\mathbf{s^3}\]
\[ +(3436
u_1^3u_2+3800u_0u_1^3+7756u_2^3u_0-3886u_1u_2^3+1225u_2^4+17059u_2^2u_0^2-5984
u_1^2u_0^2\]
\[ +15708u_2u_0^3-12232u_0u_1u_2^2+5180u_0^4-2091u_1^2u_2^2-6828u_0u_1^2
u_2+1316u_1^4-12700u_0^2u_1u_2-4312u_1u_0^3)\mathbf{s^2} \] 
\small 
$+\boldsymbol{( 384u_0^3u_1 -1792u_0^3u_2 +512u_0^2u_1^2 +1536u_0^2u_1u_2
+1920u_0u_1u_2^2 -1288u_0u_2^3 -768u_1^3u_2 }$\\ 
$\boldsymbol{-448u_0^4 -2436u_0^2u_2^2 -384u_0u_1^3 +1024u_0u_1^2u_2
-64 u_1^4 +260u_1^2u_2^2 +768u_1u_2^3 -196u_2^4)}$\normalsize
$\boldsymbol{s}.$\\ 
So our toric perturbation $\pert_{A,F^*}$ is just the 
coefficient of $s$ or $s^2$ in this polynomial, according as $\ch \K\!\neq\!2$ or $\ch \K\!=\!2$.  

Let us now examine $\pert_{A,F^*}$ itself in detail: 
Factoring with {\tt Maple}, we obtain that $\pert_{A,F^*}$ splits as  
follows:
\[ -4(u_0+u_1+u_2)(28u_0+4u_1+49u_2)(u_0-u_1+u_2)(4u_0-4u_1+u_2) 
\] 
In particular, given any factor above, the ratio of the coefficients 
of $u_i$ and $u_0$ is precisely the $i^\thth$ coordinate of some corresponding 
root of $F$. Thus the first two factors correspond precisely to the 
two isolated roots we already know.  As for the last two factors, 
note that they both give isolated points lying on the aforementioned 
line $\{-1\}\!\times\!\K$. We can then guess that this 
line should be assigned an excess intersection multiplicity of 
$2$. Of course, we might not know at the outset which of these roots is 
isolated, i.e., a zero-dimensional component of $\cZ$. However, 
since the constant term of $\cH(s)$ vanishes, assertions  
(1) of Main Theorems \ref{main:chow} and \ref{main:pert} at least tell us 
that $\cZ$ is indeed positive-dimensional. 

To distinguish the isolated roots, let us employ an algorithm from the 
proof of Corollary \ref{cor:double}: Apply Main Theorem \ref{main:fill} 
once more 
to pick $F^{**}\!=\!(1+x^3y,xy+2x^2)$. Noting that (due to their 
second equations) $F^*$ and $F^{**}$ will have no roots 
in common in $(\Ks)^2$, let us then define the {\bf double toric 
perturbation}, $\pert^{**}_A$, to be the greatest common divisor of 
$\pert_{A,F^*}$ and $\pert_{A,F^{**}}$. 

Repeating the same calculation we used for $h,h_1,h_2$, but with 
$\pert^{**}_A$ instead, we obtain new polynomials $h^{**},h^{**}_1,
h^{**}_2$. Let us compute the gcd, $g^{**}$, of $h^{**}$ and 
$h^{**}_1h^{**}_2$. It then turns out that the number of 
isolated roots of $F$ is at most $\deg h^{**}-\deg g^{**}$ 
(cf. section \ref{sub:double}). 

More explicitly, via {\tt Maple} again, we easily see that 
$h^{**}(t)\!=\!(2t+1)(2t+3)$ and $g^{**}(t)\!=\!1$. So 
the number of isolated roots in $(\Ks)^2$ is at most $2$, and the 
positive-dimensional part of $\cZ$ (the line 
$\{-1\}\times \K$) should be assigned an 
intersection multiplicity of at least $\cM(E)-2\!=\!2$. Fortuitously 
(conjecturally always), our lower bound is actually an equality.  

For completeness, we now reveal $\pert_{A,F^{**}}$ (up to 
a constant multiple): 
\[(u_0+u_1+u_2)(28u_0+4u_1+49u_2) 
(u_0-u_1+\frac{\frac{1}{\sqrt{-3}}-1}{4}u_2)
(u_0-u_1-\frac{\frac{1}{\sqrt{-3}}+1}{4}u_2).\] 
(In particular, $\pert_{A,F^{**}}$ is again the coefficient 
of $s$ in $\cH(s)$.) Note also that the last two factors 
of this toric perturbation again correspond to roots lying 
on the line $\{-1\}\!\times\!\K$.  
We thus see that varying the coefficients of our perturbation of $F$ has 
moved two of our points lying in $\cZ$. 

Note (via {\tt Maple} again) that the original GCP could have been used above, but would have resulted in a variant of 
$\pert_A$ of degree $16$ (the product of the degrees of $f_1$ and $f_2$) --- 
{\bf four times larger} than the degree of our $\pert_A$. Also, the 
old GCP is significantly larger, having {\bf 672} terms, compared to 
{\bf 110} for our above toric GCP $\cH(u;s)$. 

\subsection{Which Compactification for $\boldsymbol{\choo_A}$?}\mbox{}\\
\label{sub:stranger}
Here we show how the twisted Chow form $\choo_A$ can vanish identically for the wrong $A$, thus giving no information about 
the roots of $F$. Along the way, we will also obtain a more precise 
visualization of the toric compacta $\Pro^3_\K$, $\cT(P)$, and 
$\cT(\oP)$. We also point out that while it is customary 
to consider the roots of $F$ in $\cT(P)$ (as in 
\cite{tfulton,gkz94,toricint}), the construction of $\choo_A$ and $\pert_A$ 
necessitate the consideration of roots in $\cT(\oP)$ as well. 

To define our next example, set $n\!=\!3$, $A\!=\!\Delta\cap \Z^3$, and 
consider the $3\!\times\!3$ system $F\!=\!(a_1yz+a_2xz+a_3xy+a_4xyz, 
b_1yz+b_2xz+b_3xy+b_4xyz, c_1yz+c_2xz+c_3xy+c_4xyz)$.
Note that the mixed volume bound for this system is $1$. 
Furthermore, it is clear that $\frac{1}{xyz}F$ is a linear system in 
$\{\frac{1}{x}, \frac{1}{y},\frac{1}{z}\}$. So by Cramer's rule, we can 
express $x$, $y$, and $z$ as ratios of $3\times 3$ determinants in the 
coefficients. 

Combining this with the product formula for toric resultants \cite{chowprod} 
(and clearing denominators) we obtain that $\choo_A$ is precisely\footnote{We
also need the fact that the Pedersen-Sturmfels formula, originally
stated only over $\C$, remains true over a general algebraically
closed field (cf.\ section \ref{sub:chow}). } 
$[423][143][124]u_0+[123][143][124]u_1+[123][423][124]u_2 
+[123][423][143]u_3$  
where the {\bf bracket} $[ijk]$ \cite{introchow} is the 
$3\times 3$ subdeterminant 
\[\det\begin{bmatrix} 
a_i & a_j & a_k \\ b_i & b_j & b_k \\ c_i & c_j & c_k \end{bmatrix} \] 
of the coefficient matrix of $F$. This compactly expressed 
resultant can be thought of as a {\bf semi-mixed} Chow 
form --- a toric resultant of a system of $n+1$ polynomials 
with $k\!\ll\!n$ distinct supports.

Now consider the specialization of $F$ to $(yz+xz+2xy+3xyz,yz+xz+4xy+9xyz,\\
yz+xz+8xy+27xyz)$.
It is then easily verified that $F$ has no roots in 
$(\Ks)^3$, but $F$ does have exactly one root\footnote{If 
$\ch \K\!\in\!\{2,3\}$ then $F$ will actually have infinitely 
many roots in $\cT(P)$. So let us assume henceforth that 
$\ch \K\!\not\in\!\{2,3\}$. (It is easy to construct similar examples when 
$\ch \K\!\in\!\{2,3\}$ as well.)} in 
$\cT(P)$. Also, in our particular example, $\cT(P)\!\cong\!\Pro^3_\K$ 
and, locally (within $\Ksn$), the isomorphism is given by 
$(x,y,z)\mapsto [\frac{1}{x}\!:\!\frac{1}{y}\!:\!\frac{1}{z}\!:\!1]$. In particular, using the latter set of 
coordinates, our one root 
of $F$ in $\cT(P)$ is exactly the point $[1\!:\!-1\!:\!0\!:\!0]$.
More to the point, $\choo_A\!\equiv\!0$ for this specialization of the 
coefficients of $F$.

A simple geometric explanation for this behavior of 
$\choo_\star(\cdot)$ is that the choice of $A$ defines a toric 
variety $\cT(A)$ into which the roots of $F$ in $\cT(\oP)$ are 
projected. (The variety $\cT(A)$ is the toric variety corresponding to a 
{\bf point set} \cite{gkz94}, and is simply the image of 
$\cT(\oP)$ under the morphism $\varphi_A$.) So depending 
on our choice of $A$, the roots of $F$ in $\cT(P)$ 
may or may not correspond to roots of $F$ in $\cT(A)$ in a 
well-defined way. For instance, in our example, $F$ actually 
has {\bf infinitely} many roots in $\cT(A)$, so Main Theorem 
\ref{main:chow} tells us that $\choo_A$ {\bf must} vanish. 

So it more useful to work within $\cT(\oP)$, since the roots 
of $F$ in $\cT(P)$ and $\cT(A)$ are actually images of 
the roots of $F$ in $\cT(\oP)$. In particular, the underlying 
algebraic maps induce projections of certain faces 
of $\oP$ (corresponding to certain parts of 
$\cT(\oP)\!\setminus\!\Ksn$) onto certain faces of $P$ and 
$\conv(A)$. Figure 2 below illustrates this, 
along with where the root $[1\!:\!-1\!:\!0\!:\!0]\!\in\!\cT(P)$ of 
$F$ ``goes'' within these various compacta. For instance, note that 
$\oP$ is a cuboctahedron, and $\fii_A$ is constant on the portions of 
$\cT(\oP)\!\setminus\!\Ksn$ corresponding to the triangular faces 
with inner normals $-\hat{e}_1$, $-\hat{e}_2$, $-\hat{e}_3$, and 
$(1,1,1)$.  

\begin{figure}[h]
\begin{center}
 \epsfig{height=3.4in,file=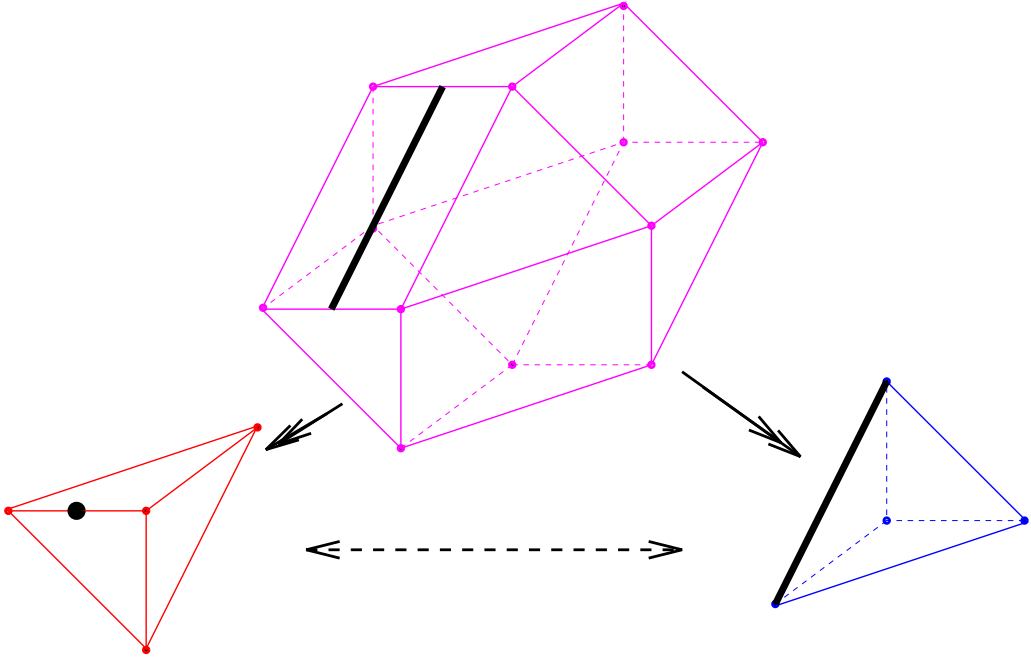}
\end{center} 
\caption{One root in the lower left toric compactification ($\cT(P)$) 
becomes infinitely many roots in the other two compactifications 
($\cT(\oP)$ and $\cT(A)$).} 
\end{figure}

Algebraically, we have the following maps: 

\begin{picture}(140,60)(-150,-5)
        \put(55,40){$\cT(\oP)$} 
        \put(57,42){\begin{rotate}{-135}
                     \begin{huge} 
                      $\twoheadrightarrow$
                     \end{huge}
                    \end{rotate}} 
        \put(38,28){\begin{scriptsize}
                     $\pi$
                    \end{scriptsize}}
        \put(74,33){\begin{rotate}{-45}
                     \begin{huge}
                      $\twoheadrightarrow$
                     \end{huge} 
                    \end{rotate}} 
        \put(88,28){\begin{scriptsize}
                     $\varphi_{_A}$
                    \end{scriptsize}}
\put(-17,0){$\Pro^3_\K\cong\cT(P)$}
             \put(39,0){\begin{Large}
                         $\dashleftarrow\dashrightarrow$
                        \end{Large}}
                           \put(100,0){$\cT(A)\hookrightarrow\Pro^3_\K$} 
        \put(62,-4){\begin{scriptsize}
                     $\phi$
                    \end{scriptsize}}
\end{picture} 

\noindent
where $\pi$ is the natural projection between compatible 
toric compacta (cf.\ section \ref{sec:chowtor}), and $\phi$ is the rational 
map (defined just on $\Ksn$) from $\cT(P)$ to 
$\cT(A)$ obtained from $x\mapsto [x^a \; | \; a\!\in\!A]$. 
In the case at hand, the latter map is simply the 
identity map between the two corresponding naturally embedded copies of 
$(\Ks)^3$. 

To remedy the preceding trivial $\choo_A$, we can instead use 
$\choo_{A'}(u)$ with $A'\!:=\!\{(0,1,1),\\ 
(1,0,1),(1,1,0),(1,1,1)\}$. 
(This choice is motivated by trying to pick an $A'$ 
which is compatible with $P$ (cf.\ section \ref{sec:chowtor}).) In particular, 
when the coefficients of $F$ are unspecialized, 
$\choo_{A'}(u)\!=\!\det\!$\begin{scriptsize} 
$\begin{bmatrix} 
  a_1 & a_2 & a_3 & a_4\\ 
  b_1 & b_2 & b_3 & b_4\\
  c_1 & c_2 & c_3 & c_4\\
  u_{(0,1,1)} & u_{(1,0,1)} & u_{(1,1,0)} & u_{(1,1,1)} 
\end{bmatrix}$\end{scriptsize}. 
So under our last specialization, this becomes 
$12u_{(1,0,1)}-12u_{(0,1,1)}$. Note that we now recover 
our root $[1\!:\!-1\!:\!0\!:\!0]$ from the coordinates of our new twisted 
Chow form. For example, the ratio of the $x$-coordinate to 
the $y$-coordinate is just 
$\frac{x}{y}\!=\!\frac{x^1y^0z^1}{x^0y^1z^1}\!=\!\frac{12}{-12}\!=\!-1$. 

Alternatively, we can simply use $\pert_A$ and forget about 
cleverly chosen $A'$. For example, by Main Theorem 3 (and theorem 
\ref{thm:me}), we can simply take $F^*\!=\!(yz+xyz,xz+xyz,xy+xyz)$. After an 
application of {\tt Maple}, we then obtain that $\pert_{\Delta\cap\Z^3}$ is 
exactly $5u_1+21u_2$. In particular, while the point 
$[0\!:\!5\!:\!21\!:\!0]\!\in\!\cT(A)$ does not correspond (in any 
obvious way) to a root of $F$ in $\cT(P)$, it {\bf is} the image of a bona 
fide root of $F$ in $\cT(\oP)$ under the morphism $\varphi_A$.  

In closing, we emphasize that in practice we would never 
actually compute the full monomial expansions of 
$\choo_A(u)$, $\pert_A(u)$, or $\cH(u;s)$ --- we would 
instead recover the roots of $F$ (or evaluate monomials thereof) via 
rapid and sophisticated interpolation 
techniques, e.g., \cite{pspace,ckl89,gcp,black}. In particular, this 
is the approach of Main Theorem 1, and our calculations can be sped up 
tremendously with suitably optimized code. 

\subsection{The ``Dense'' Case}\mbox{}\\
\label{sub:dense}
Our last example illustrates a simple fundamental case.

Suppose $E$ is the $n$-tuple $(d_1\Delta\!\cap\!\Zn,\ldots,
d_n\Delta\!\cap\!\Zn)$ where $d_i\!\in\!\N$ for all $i$. (So 
we are now considering the family of all $n\!\times\!n$ polynomial systems 
where $f_i$ has total degree $\leq\!d_i$ for all $i$.) It is then easily 
verified that the system $F^*\!=\!(x^{d_1}_1,\ldots,x^{d_n}_n)$ (with 
support contained in $E$) has only finitely many roots in $\cT(P)$. 
Indeed, in this case, $\cT(P)\!\cong\!\Pn$ and there 
is exactly one root (of multiplicity $\prod d_i$) at the 
origin $\bO$. Note also that our current setting is sufficiently 
simple that we could find a suitable $F^*$ with just $n$ terms, 
without the need for an irreducible fill. 

\begin{rem} 
Letting $d_\Pi\!:=\!\prod d_i$ and $d_\Sigma\!:=\!\sum^n_{i=1} d_i$, it 
is easily checked that $\cM(E)\!=\!d_\Pi$ and 
$R(\oE)\!=\!S(\oE)\!=\!$\begin{tiny}$\begin{pmatrix} d_\Sigma+1\\ n 
\end{pmatrix}$\end{tiny} in the dense case. The last 
equality follows from Macaulay's 19$^\thth$ century construction 
of the multivariate resultant \cite{cannyphd}. So in the dense case, 
our complexity bound from Main Theorem 1 thus specializes to the best 
bounds from \cite{pspace,ckl89,gcp}, once $\ch \K\!=\!0$ and 
randomization is allowed 
(cf.\ Corollary \ref{cor:fast} of section \ref{sec:complex}). 
\end{rem} 

Letting $A\!=\!\Delta\!\cap\!\Zn$, we then see that our polynomial $\cH(u;s)$ 
is simply the original GCP \cite{gcp}, but extended to a general algebraically 
closed field. In particular, our $F-sF^*$ 
is the polynomial system $(f_1-sx^{d_1},\ldots,f_n-sx^{d_n})$.  
(Note also that if we set $d_1\!=\cdots=\!d_n\!=\!1$ 
then we recover the usual characteristic polynomial of a matrix). 
Finally, note that $\cT(A)\!\cong\!\cT(\oP)\!\cong\!\cT(P)\!\cong\!\Pn$ 
and the map $\fii_A$ is the identity. So by 
considering the zero set of $F$ in $\cT(\oP)$, in this ``dense'' case, we are 
just considering the zero set of $F$ in $\Pn$ in the usual way via 
homogenizations. Thus by Main Theorem \ref{main:pert}, Canny's original GCP 
indeed finds a point in every irreducible component of $\cZ$ in $\Pn$, as 
conjectured in 1990. Of course, the advantage of the toric GCP 
is that we can do the same with greater efficiency for sparse 
systems with small $\cM(E)$. 

\section{Filling}
\label{sec:fill}
Here we briefly recount filling and some
related concepts. Some of the material below is covered at greater length
in \cite{convexapp}. The results below form the basis for our combinatorial 
approach to perturbing degenerate polynomial systems.  

Let $\Sn\!\subset\!\Rn$ denote the unit $(n\!-\!1)$-sphere centered at the
origin. For any compact $B\!\subset\!\Rn$ and any $w\!\in\!\Rn$, define
$B^w$ to be the set of $x\!\in\!B$ where the inner-product $x\!\cdot\! w$ is
minimized. (Thus $B^w$ is the intersection of $B$ with its supporting
hyperplane in the direction $w$.) We then define
$E^w\!:=\!(E^w_1,\ldots,E^w_n)$ and
$D\!\cap\!E^w\!:=\!(D_1 \cap E^w_1,\ldots,D_n \cap E^w_n)$.

Recall that the {\bf dimension} of any $B\!\subseteq\!\Rn$, $\dim B$, is
the dimension of the smallest subspace of $\Rn$ containing a translate of
$B$. The following definition is fundamental to our development.
\begin{dfn}
\label{dfn:ess} 
Suppose $C\!:=\!(C_1,\ldots,C_n)$ is an
$n$-tuple of polytopes in $\Rn$ {\bf or} an $n$-tuple of finite subsets of
$\Rn$. We will allow any $C_i$ to be empty and say that a nonempty subset
$J\!\subseteq\![n]$ is {\bf essential} for $C$ (or $C$ {\bf has
essential subset $J$}) $\Longleftrightarrow$ (0) $\supp(C)\!\supseteq\!J$,
(1) $\dim(\sum_{j\in J} C_j)=\#J-1$, and (2) $\dim(\sum_{j\in J'} C_j)\geq
\#J'$ for all nonempty {\bf proper} $J'\!\subsetneqq\! J$.
\end{dfn}

Equivalently, $J$ is essential for $C \Longleftrightarrow$ the
$\#J$-dimensional mixed volume of $(C_j \; | \; j\!\in\!J)$ is $0$ and no
smaller subset of $J$ has this property. Figure \ref{fig:myfirst} below shows 
some simple examples of essential subsets for $C$, for various $C$ in the case
$n\!=\!2$. 

\begin{figure}[h]
\begin{picture}(410,70)(-40,-5)

\put(-80,-50){
\begin{picture}(133,100)(0,0)
\put(51,80){\circle*{3}} \put(45,85){$C_1$} 
\put(72,80){\circle*{3}} \put(69,85){$C_2$} 
\put(42,50){$\{1\},\{2\}$}
\end{picture}}

\put(24,-50){
\begin{picture}(133,100)(0,0)
\put(51,80){\circle*{3}} \put(45,85){$C_1$}
\put(62,70){\line(1,1){30}}
\put(62,70){\circle*{3}} \put(92,100){\circle*{3}} \put(79,77){$C_2$}
\put(60,50){$\{1\}$}
\end{picture}}

\put(145,-50){
\begin{picture}(133,100)(0,0)
\put(36,70){\line(1,1){30}}
\put(36,70){\circle*{3}} \put(66,100){\circle*{3}} \put(36,85){$C_1$}
\put(62,70){\line(1,1){30}}
\put(62,70){\circle*{3}} \put(92,100){\circle*{3}} \put(79,77){$C_2$}
\put(45,50){$\{1,2\}$}
\end{picture}}

\put(276,-50){
\put(26,100){\line(1,-1){30}}
\put(26,100){\circle*{3}} \put(56,70){\circle*{3}} \put(43,85){$C_1$}
\put(72,70){\line(0,1){30}}
\put(72,70){\circle*{3}} \put(72,100){\circle*{3}} \put(74,85){$C_2$}
\begin{picture}(133,100)(0,0)
\put(45,50){None}
\end{picture}}

\end{picture}

\caption{The essential subsets for 4 different pairs of plane polygons. 
(The segments in the third pair are meant to be parallel.) } 
\label{fig:myfirst}

\end{figure}
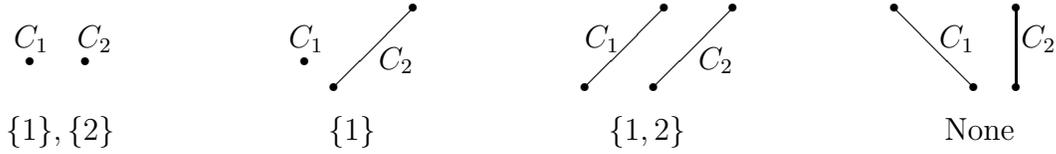

A basic fact about mixed volumes is that $\cM(E)\!=\!0 \Longleftrightarrow
E$ has an essential subset, whenever $\supp(E)\!=\![n]$.
However, there is an even deeper connection between filling and 
essentiality:
\begin{thm} \cite[sec.\ 2.5]{convexapp}
\label{thm:me}
Suppose $D$ and $E$ are $n$-tuples of finite subsets of $\Zn$ such that
$\cM(E)\!>\!0$. Then $D$ fills $E
\Longleftrightarrow$ for all $w\!\in\!\Sn$, $\supp(D\cap E^w)$
contains a subset essential for $E^w$. \qed 
\end{thm}
\begin{rem}
One certainly need not check infinitely many $w$. In fact, 
we need only check one $w$ (just pick any inner normal) 
for each face of the polytope $P\!=\!\sum^n_{i=1} \conv(E_i)$.
\end{rem}

We also recall the following important observation. 
\begin{lemma}
\label{lemma:ed} 
Let $m\!:=\!\sum^n_{i=1}\#E_i$. Then there is an  
absolute constant $c\!>\!0$ such that we can 
decide whether $\cM(E)\!>\!0$ within $\cO(m^c)$ arithmetic 
steps over $\Q$. Furthermore, if $\cM(E)\!=\!0$, 
then we can find points $p_1,\ldots,p_n\!\in\!\Non$ such 
that $\cM(\{p_1\}\cup E_1,\ldots,\{p_n\}\cup E_n)\!>\!0$ within 
the same asymptotic complexity bound. \qed 
\end{lemma} 
\noindent 
The first portion appears in \cite{gk94}, while the 
second portion is an elementary corollary. 

Oddly enough, filling seems to have originated from an 
algebraic problem: genericity conditions for counting the roots of sparse 
polynomial systems. This aspect is explored much further in
\cite{convexapp,rojaswang,toricint}. We also emphasize that constructing a 
fill need only be done {\bf once} for a given family of problems, provided
$E$ remains fixed. The situation where the monomial term structure of a
polynomial system remains fixed once and for all, and the coefficients may vary
many thousands of times, actually occurs frequently in many practical 
contexts such as robot control or computational geometry.

To conclude our background, we will need the following lemma characterizing 
irreducible fills. 
\begin{lemma}
\label{lemma:gen}
Following the preceding notation, assume $\cM(D)\!>\!0$. Then $D$ is 
irreducible $\Longleftrightarrow$ for any $v$ lying in some $D_i$, there exists a $w\!\in\!\Qns$ such that $D^w_i\!=\!\{v\}$ and 
$\cM(D^w_1,\ldots,D^w_{i-1},D^w_{i+1},\ldots,D^w_n)\!>\!0$. 
\end{lemma}
\noindent
{\bf Proof:} First note that the mixed volume condition 
above is equivalent to $\{i\}$ being the unique essential 
subset of $D^w$. This follows immediately from definition \ref{dfn:ess} 
and, say, the development of \cite{buza}. 

The ``$\Longleftarrow$'' direction then follows 
almost immediately from theorem \ref{thm:me}: If the mixed volume condition 
holds, then the removal of any point from $D$ would indeed violate the filling 
condition from theorem \ref{thm:me}. So the removal of any point from 
$D$ would make $\cM(D)$ decrease. The converse implication follows almost as 
easily: 

Suppose, to derive a contradiction, that $D$ is irrreducible 
but there is some $v$ in some $D_i$ satisfying the following 
property: For all $w\!\in\!\Qns$, $\#D^w_i\!\geq\!2$ or 
$\cM(D^w_1,\ldots,D^w_{i-1},D^w_{i+1},\ldots,\\ D^w_n)\!=\!0$. 
Let us then consider the $n$-tuple $D'\!:=\!(D_1,\ldots,
D_{i-1},D_i\!\setminus\!\{v\},D_{i+1},\ldots,D_n)$. 
Then by theorem \ref{thm:me} once again, $D'$ fills 
$D$. But this contradicts the irreducibility of $D$, so 
we are done. \qed 

\section{Toric Geometry and the Proofs of Our Main Theorems} 
\label{sec:chowtor} 
Our notation is a slight variation of that used in 
\cite{tfulton}, and is described at greater length in \cite{toricint}.  
However, we will briefly review a few important facts and definitions. 

The {\bf (inner) normal fan} of a polytope 
$Q\!\subset\!\Rn$, $\fan(Q)$, is simply the collection of cones of 
inner normals of faces of $Q$ \cite{gkz94}. (For instance, the 
inner normal fan of the standard unit square in the plane consists 
of nine cones: the four quadrants, the four nonnegative coordinate 
rays, and the origin.) We will assume the reader to be familiar with 
the construction of a toric variety from a fan, a polytope, or a finite point 
set \cite{tfulton,gkz94}. 
\begin{ex} When $A\!=\!\Delta\cap\Zn$, it is easy to derive from 
scratch that $\cT(A)$ is just the projective 
space $\Pn$. More generally, if $\conv(A)$ is a product of 
simplices, then $\cT(A)$ is a {\bf product} of {\bf twisted} projective 
spaces \cite{tfulton} --- hence our appelation for $\choo_A$. Note also 
that the coefficients of $\choo_A$ are multisymmetric functions   
of $\{\fii_A(\zeta)\}_\zeta$ as $\zeta$ ranges over the roots of $F$ 
in $\cT(\oP)$. 
\end{ex} 

Let us now list our cast of main characters:
\begin{dfn} \cite{tfulton,toricint} 
\label{dfn:toric}
Given any $w\!\in\!\Rn$, we will use the following notation:
\begin{itemize}
\item[$T=$]{The algebraic torus $\Ksn$}
\item[$Q^w=$]{The face of $Q$ with inner normal $w$} 
\item[$\sigma_w=$]{The closure of the cone generated by the inner normals of 
$Q^w$} 
\item[$\sigma^\vee_w$=]{The {\bf dual} (or {\bf angle}) cone 
$\{w'\!\in\!\Rn \; | \; w'\!\cdot\!y\!\geq\!0 \text{ for all } 
y\!\in\!\sigma_w\}$.} 
\item[$U_w=$]{The affine chart of $\cT(Q)$ corresponding to all 
semigroup homorphisms\footnote{Note that the domain and range spaces are  
respectively semigroups under the natural operations of vector 
addition and field multiplication.} $\sigma^\vee_w\cap\Zn \longrightarrow 
\K$. }  
\item[$x_w=$]{The point in $U_w$ corresponding to the semigroup
homomorphism $\sigma^\vee_w\cap\Zn \longrightarrow \{0,1\}$ mapping $p
\mapsto \delta_{w\cdot p,0}$, where $\delta_{ij}$ denotes the Kronecker
delta}  
\item[$O_w=$]{The $T$-orbit of $x_w=$ The $T$-orbit corresponding to 
the relative interior of $Q^w$} 
\item[$\cE_Q(Q')=$]{The $T$-invariant Weil divisor of $\cT(Q)$ corresponding to 
a polytope $Q'$.} 
\item[$\divisor(f)=$]{The Weil divisor of $\cT(Q)$ defined by 
a rational function $f$ on $\Ksn$} 
\item[$\cD_Q(f,Q')=$]{$\divisor(f)+\cE_Q(Q')=$ The toric effective divisor of 
$\cT(Q)$ corresponding to $(f,Q')$} 
\item[$\cD_Q(F,\cP)=$]{The (nonnegative) cycle in the Chow ring of 
$\cT(Q)$ defined by $\bigcap^k_{i=1} \cD_Q(f_i,P_i)$, whenever 
$\cP\!=\!(P_1,\ldots,P_k)$} 
\end{itemize}
\end{dfn} 
We say that $P$ is {\bf compatible} with $Q$ iff 
every cone of $\fan(Q)$ is a union of cones of $\fan(P)$ 
\cite{khocompat,tfulton,
toricint}. (So $P$ compatible with $Q \Longrightarrow P$ has 
at least as many facets as $Q$.) 
\begin{figure}[h]
\begin{center}
\epsfig{height=1.5in,file=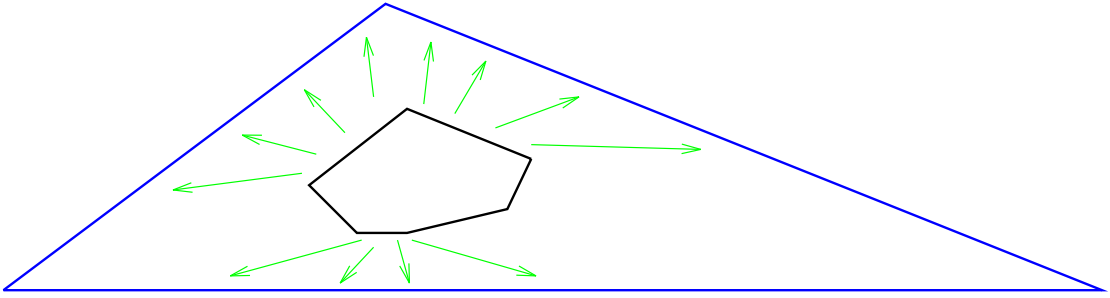}
\end{center}
\caption{The inner polytope is compatible with 
the outer polytope. Also, the corresponding ``outer'' 
toric variety can be obtained as a deformation (or image 
under a proper morphism) of the ``inner'' toric variety. }
\end{figure}

Finally, whenever $F$ is a $k\times n$ polynomial system with support contained in $E$, we will define the {\bf zero set\footnote{ 
When necessary, we will also use the underlying scheme structure.} of 
$\boldsymbol{F}$ in $\boldsymbol{\cT(Q)}$} to be the toric cycle 
$\cD_Q(F,\cP)$, where $\cP\!:=\!(\conv(E_1),\ldots,\conv(E_k))$. 
{\bf Toric infinity} is then defined relative to $Q$: 
it is simply the set $\cT(Q)\!\setminus\!\Ksn$. 
\begin{ex} {\bf (Zero Sets in $\Pn$)}  
\label{ex:pn} 
Suppose $Q\!=\!\alpha+\beta\Delta$, for any $\alpha\!\in\!\Qn$ 
and any rational $\beta\!>\!0$. Then $\cT(Q)\!\cong\!\Pn$ 
canonically. As for explicitly defining the zero set of $F$ in 
$\cT(Q)$, we can do the following: (1) Define vectors 
$p_1,\ldots,p_n\!\in\!\Zn$ 
such that for all $i$, $x^{p_i}f_i\!\in\!\K[x]$ is {\bf not} 
divisible by any $x_j$, (2) define 
$\tilde{f}_i(x)\!:=\!x^{d_i}_\infty x^{p_i} 
f(\frac{x_1}{x_\infty},\ldots,\frac{x_n}{x_\infty})$ for all 
$i$, where $d_i$ is the total degree of $x^{p_i}f_i$. Then 
$[z_1:\cdots :z_n:z_\infty]\!\in\!\Pn$ is a root 
of $F$ iff $\tilde{f}_1(z)\!=\cdots=\!\tilde{f}_n(z)\!=\!0$.  
In particular, note that this toric definition differs from the classical 
definition of ``zero set of $F$ in $\Pn$,'' due to the 
extra step (1). For instance, our toric definition might omit 
some affine roots, for certain $E$ and $F$. However, note that step (1) is 
unnecessary when $\bO\!\in\!E_i$ for all $i$. Furthermore, by 
\cite[sec.\ 6.1]{toricint}, the zero scheme of $F$ in $\Kn$ embeds naturally in $\cD_P(F,\cP)$ (and $\cD_\oP(F,\cP)$) when we replace $\oE$ by $\bO\cup\oE$. 
Hence the introduction of $\bO\cup E$ in (and $A\!=\!\Delta\cap\Zn$ in 
the proof of) Main Theorem 1. 
\end{ex} 

The following result will provide some necessary geometric intuition 
for specializing resultants. The lemma immediately following then gives 
a more explicit algebraic analogy between the faces of $Q$ and the affine 
charts of $\cT(Q)$.  
\begin{van} \cite{resvan} 
Suppose $\oF$ is an $(n+1)\times n$ polynomial system (over $\K$) with support 
contained in $\oE$. Then, provided 
$\cM(E_1,\ldots,E_{i-1},E_{i+1},\ldots,E_{n+1})\!>\!0$ 
for some $i\!\in\![n+1]$, $\res_{\bar{E}}(\oF)\!=\!0 
\Longleftrightarrow \cD_\oP(F,\cP)\!\neq\!\emptyset$,
where $\cP\!:=\!(\conv(E_1),\ldots,\conv(E_{n+1}))$ and 
$\oP\!=\!\sum^{n+1}_{i=1} \conv(E_i)$. \qed 
\end{van}
\begin{lemma} \cite[sec.\ 4.2--5.1]{toricint}
\label{lemma:once}
Suppose $F$ is a $k\times n$ polynomial system over $\K$ with support 
contained in a $k$-tuple of integral polytopes $\cP\!:=\!(P_1,\ldots,P_k)$ in 
$\Rn$. Assume further that $Q$ is a rational polytope in $\Rn$. 
Then the defining ideal in $\K[x^a \; | \; a\!\in\!\sigma^\vee_w\cap \Zn]$ of $U_w\cap \cD_Q(F,\cP)$ is $\langle x^{b_i}f_i \; | \;$ for all 
$i\!\in\![k]$ and $b_i\!\in\!\Zn$ such that $b_i+P_i\subseteq 
\sigma^\vee_w \rangle$. \qed
\end{lemma} 

Lifting (or projecting) from one toric variety to another is an 
important fundamental ideal we will also use. The following lemma 
follows directly from the development of \cite{tfulton}.
\begin{lemma}
\label{lemma:lift}
Suppose $Q\!\subset\!\Rn$ is an $n$-dimensional rational polytope, 
and $B$ is either a nonempty finite subset of $\Zn$ {\bf or} a 
rational polytope in $\Rn$. Assume further that $Q$ is compatible 
with $\conv(B)$. Then 
there is a natural (surjective) proper morphism $\pi : \cT(Q) 
\twoheadrightarrow \cT(B)$. In particular, 
$\pi(\cD_Q(F,\cP))\!=\!\cD_B(F,\cP)$, where the latter 
cycle is the image of $\cD_{\conv(B)}(F,\cP)$ under the natural 
proper morphism from $\cT(\conv(B))$ to $\cT(B)$. 
Furthermore, $\pi(O_w)\!=\!O_w$, where the corresponding $T$-orbits are 
considered in their respective domains, and 
\mbox{$\pi|_{\Ksn}\!=\!\mathrm{id}$. \qed } 
\end{lemma} 
\begin{rem} 
\label{rem:mult} 
Following the notation of Main Theorem 1, it easily follows that if 
$A\!=\!\Delta\cap\Zn$ then the multiplicity of any root of $F$ in $\Ksn$ 
is preserved under the map $\fii_A$. If $A\!=\!dQ\cap\Zn$ 
for some rational polytope compatible with $P$, and $d\!\in\!\N$ is 
sufficiently large, then the same will be 
true of any root of $F$ in $\cT(\oP)$ \cite{tfulton}. In general, 
thanks to the functoriality of Chow forms \cite{introchow}, 
$\choo_A$ is precisely the Chow form of the 
subscheme $\fii_A(\cZ)$ of $\Pa$. 
\end{rem} 

Another immediate corollary of our last lemma is the following 
result on the meaning of the projective coordinates $[\zeta^a \; 
| \; a\!\in\!A]$.
\begin{cor}
\label{cor:parts} 
Following the notation of Main Theorem 4, let 
$\zeta\!\in\!\cT(\oP)$ be an isolated root of $F$ and fix a 
vertex of $\mathrm{v}\!\in\!\conv(A)$ with inner normal $w$. 
Then $\fii_A(\zeta)$ lies in the affine chart 
$U_w$ of $\cT(A) \Longleftrightarrow$ the 
coefficient of $u_{\mathrm{v}}$ in the corresponding factor 
of $\pert_A$ is {\bf non}zero. \qed 
\end{cor}
\begin{ex}
Suppose we take $A\!=\!\Delta\cap\Zn$ as usual. Then 
$\cT(A)\!\cong\!\Pn$ canonically, and there are exactly 
$n+1$ affine charts of $\cT(A)$ corresponding to vertices. These 
charts are respectively isomorphic to $\Pn$ minus the hyperplane at infinity, 
and $\Pn\!\setminus\!\{x_i\!=\!0\}$ as $i$ runs through $[n]$. For example, 
given a factor of $\pert_A$ such as $u_0+u_3$, we know that it corresponds 
to a root image $\fii_A(\zeta)$ which lies in 
two of these affine charts and outside of $n\!-\!2$ others, i.e., 
$\fii_A(\zeta)\!=\![0\!:\!0\!:\!1\!:\!0\!:\cdots:\!0\!:\!1]$ lies on the 
$x_3$-axis. Similarly, if all the coordinates of $\fii_A(\zeta)$ are 
{\bf non}zero, then $\zeta,\fii_A(\zeta)\!\in\!\Ksn$. 
\end{ex} 

Finally, we will need a version of the fundamental 
fact that $F$ generically has exactly $\cM(E)$ roots 
in $\Ksn$. The case $\K\!=\!\C$ first appeared in 
\cite{bernie}, and the general case is an immediate 
corollary of \cite[Main Theorem 2]{toricint}. 
\begin{lemma}
\label{lemma:bernie} 
Let $\cC_E$ be the vector of coefficients 
of $F$ and define $\#E\!:=\!\sum^n_{i=1} \#E_i$. Then 
there is an algebraic hypersurface 
$\Sigma_E\!\subset\!\K^{\#E}$ such that 
$\cC\!\in\!\K^{\#E}\!\setminus\!\Sigma_E \Longrightarrow 
F$ has {\bf no} roots in $\cT(P)\!\setminus\!\Ksn$. 
Moreoever, the latter assertion implies that $F$ has 
exactly $\cM(E)$ roots, counting multiplicities, in 
$\Ksn$. \qed 
\end{lemma} 

With all our technical background complete, we can now prove 
our main theorems. 

\subsection{Polynomial Algebra and the First Half of Main Theorem 1} 
\label{sub:solve} \mbox{}\\
Our proof of assertions (0)--(2) of Main Theorem 1 will rely on two main constructions: the toric 
perturbation $\pert_{\Delta\cap\Zn}$ and an 
extension of Canny's constructive version \cite{pspace} of the 
primitive element theorem. 
We thus emphasize that while $\choo_A$ and $\pert_A$ permit 
one to reduce polynomial system solving to multivariate 
factorization, we will {\bf not} use factoring to 
build $h$ and $h_1,\ldots,h_n$. 

Algebraically, the idea is as follows: Our techniques 
allow us to find a set of points $Z'\!\subset\!\Ksn$  
intersecting every irreducible component of the zero set of 
$F$ in $\Ksn$.  
Consider the field extension $L\!:=\!K(Z')$, obtained by 
adjoining all the coordinates of all the points of 
$Z'$. Then $L$ is a finite extension of $K$, and by the 
primitive element theorem \cite{vdv}, $L\!=\!K(\theta)$ 
for some $\theta\!\in\!L$. Furthermore, by the same theorem, we 
should be able to recover the coordinates of every point in $Z'$ in 
terms of rational functions (with coefficients in $K$) of $\theta$. 
Since $K(\theta)\!\cong\!K[t]/h(t)$ when $h$ is the minimal 
polynomial of $\theta$ over $K$, we can further simplify 
the preceding rational representation to one in terms of polynomials in 
$\theta$ with coefficients in $K$.  Our algorithm for Main Theorem 
\ref{main:big} will explicitly construct this encoding for us. 

To describe our algorithm, we will first need a bit of 
subresultant theory: For any univariate polynomials 
$f(t)\!=\!\alpha_0+\alpha_1t+\cdots+\alpha_{d_1}t^{d_1}$ and 
$g(t)\!=\!\beta_0+\beta_1t+\cdots+\beta_{d_2}t^{d_2}$, consider the following 
$(d_1+d_2-2)\times (d_1+d_2-1)$ matrix 
\begin{scriptsize} 
\[ 
\begin{bmatrix} 
\beta_0 & \cdots & \beta_{d_2} & 0   & \cdots & 0 & 0 \\
0      & \beta_0 & \cdots & \beta_{d_2} & 0 & \cdots & 0\\
\vdots & \ddots & \ddots &  & \ddots & \ddots & \vdots \\
0      & \cdots & 0 & \beta_0 & \cdots & \beta_{d_2} & 0 \\
0      & 0  & \cdots & 0 & \beta_0 & \cdots & \beta_{d_2} \\ 
\alpha_0 & \cdots & \alpha_{d_1} & 0 & \cdots & 0 & 0 \\
0     & \alpha_0 & \cdots & \alpha_{d_1} & 0 & \cdots & 0\\
\vdots &  \ddots  & \ddots &    &  \ddots  & \ddots  & \vdots   \\
0      & \cdots & 0 & \alpha_0 & \cdots & \alpha_{d_1} & 0 \\
0      & 0  & \cdots & 0 & \alpha_0 & \cdots & \alpha_{d_1} 
\end{bmatrix}
\] 
\end{scriptsize}
\hspace{-\sh}with $d_1\!-\!1$ ``$\beta$ rows'' and $d_2\!-\!1$ ``$\alpha$ 
rows.'' Let $M^1_1$ (resp.\ $M^1_0$) be the submatrix obtained by 
deleting the last (resp.\ second to last) column, 
and let $\cR_i(f,g)\!:=\!\det(M^1_i)$ for $i\!\in\!\{0,1\}$. Finally, 
define the {\bf first subresultant} of $f$ and $g$ to be 
$\cR_0(f,g)+\cR_1(f,g)t$. It is then a classical 
fact that if $\gcd(f,g)\!=\!a+bt$ with $b\!\neq\!0$, then 
$\frac{a}{b}\!=\!\frac{\cR_1(f,g)}{\cR_0(f,g)}$ \cite{lalo}. We will make 
heavy use of this fact in our proof. 

Recall also the following algorithmic facts about polynomials 
over any field \cite{binipan}: 
\begin{enumerate} 
\renewcommand{\theenumi}{\alph{enumi}}  
\item{Given the values of a univariate polynomial of degree $d$ at 
$d+1$ distinct points, the coefficients of the polynomial 
can be recovered within $\cO^*(d)$ field operations. }
\item{The gcd of two univariate polynomials of 
degree $\cO(d)$ can be found within $\cO^*(d)$ field 
operations.} 
\item{The coefficients of the square-free part of a univariate polynomial (of 
degree $d$) can be found within $\cO^*(d)$ field operations. } 
\item{The subresultant of two univariate polynomials of degree 
$\cO(d)$ can be computed within $\cO^*(d)$ additions and 
multiplications.}  
\end{enumerate} 

We now proceed with our proof of the first half of Main Theorem \ref{main:big}.

\noindent
{\bf Proof of Assertions (0)--(2):}
To simplify matters slightly, we will first derive a Las Vegas 
version of our algorithm for Main Theorem \ref{main:big}. The announced time 
bound will then follow from a simple derandomization.
The construction of $h,h_1,\ldots,h_n$ will follow 
from evaluating $\pert_A(u)$ at various specializations 
of $u$, thus reducing to $\cO(n)$ univariate polynomial interpolation 
and gcd problems. In particular, our algorithm can be outlined as follows:
\begin{itemize}
\item[{\bf Step 0}]{Set $A\!=\!\Delta\cap\Zn$ and fix {\bf generic} values in 
$\K$ for $u_1,\ldots,u_n$.}  
\item[{\bf Step 1}]{Define $h\!\in\!\K[t]$ to be $\pert_A(t,u_1,\ldots,u_n)$.} 
\item[{\bf Step 2}]{If $n\!=\!1$, set $h_1(\theta)\!:=\!\theta$ and stop. 
Otherwise, for all $i\!\in\![n]$, let $q^-_i(t)$ be the square-free part of 
$\pert_A(t,u_1,\ldots,u_{i-1},u_i-1, u_{i+1},\ldots,u_n)$.}
\item[{\bf Step 3}]{Let $\alpha$ satisfy either $\alpha\!=\!1$ or 
$\alpha(\alpha+1)\!=\!1$ according as $\ch \K\!\neq\!2$ or 
$\ch\K\!=\!2$. Then define $q^\star_i(t)$ to be the square-free part of 
$\pert_A(t,u_1,\ldots,u_{i-1},u_i+\alpha,u_{i+1},\ldots,u_n)$ for all 
$i\!\in\![n]$.}  
\item[{\bf Step 4}]{For all $i\!\in\![n]$ and $j\!\in\!\{0,1\}$, let 
$r_{i,j}(\theta)$ be the reduction of $\cR_j(q^-_i(t),
q^\star_i((\alpha+1)\theta-\alpha t))$ modulo $h(\theta)$. }  
\item[{\bf Step 5}]{For all $i\!\in\![n]$, define 
$h_i(\theta)$ to be the reduction of 
$-\theta-\frac{r_{i,1}(\theta)}{r_{i,0}(\theta)}$ modulo $h(\theta)$.} 
\end{itemize} 
\noindent
Note that assertion (0) thus follows immediately from 
Steps 1 and 4, thanks to the beginning of Main Theorem \ref{main:pert}. 
Let us now verify the correctness of our algorithm, clarifying the 
genericity assumption of Step 0 along the way. 

Using Main Theorem 4 once more, we know that the factors of $\pert_A$ define 
for us a set of points $Z\!=\!\{\zeta^{(j)}\}_{j\in [N]}$, with 
$N\!\leq\!\cM(E)$, such that $Z$ intersects every irreducible component of 
$\fii_A(\cZ)$. In particular, we see that the roots of $h$ are 
exactly $\{\theta^{(j)}\}_{j\in N'}$, where 
$\theta^{(j)}\!:=\!-\sum^n_{i=1}\zeta^{(j)}_iu_i$, 
$Z'\!:=\!\{\zeta^{(j)}\}_{j\in N'}\!=\!Z\cap\Kn$,  
and $\zeta^{(j)}\!=\!(\zeta^{(j)}_1,\ldots,\zeta^{(j)}_n)$ for all 
$j\!\in\!N'$. Furthermore, it is easy to check that for all but finitely many $[u_1:\cdots:u_n]$, $j\!\neq\!j' \Longrightarrow 
\theta^{(j)}\!\neq\!\theta^{(j')}$. (In which case, via 
remark \ref{rem:mult} in section \ref{sec:chowtor}, the multiplicity
of any {\bf isolated} root $\zeta^{(j)}\!\in\!\Ksn$ of $F$ is 
exactly the multiplicity of the root $\theta^{(j)}$ of $h$.) 
Similarly, for any $i\!\in\![n]$, $j\!\neq\!j' \Longrightarrow 
\theta^{(j)}+\zeta^{(j)}_i\!\neq\!\theta^{(j')}+\zeta^{(j')}_i$ 
and 
$\theta^{(j)}-\alpha\zeta^{(j)}_i\!\neq\!\theta^{(j')}-\alpha\zeta^{(j')}_i$,  
for all but finitely many $[u_1:\cdots:u_n]$. The avoidance of 
these $1+2n$ finite sets of $[u_1:\cdots:u_n]$ is 
precisely our genericity condition for Step 0. Furthermore, by checking 
square-free parts, we can check our genericity condition 
with negligible overhead (via fact (c)).  

Now note that if $\theta\!=\!\theta^{(j)}$ for some $j$, then 
for all $i\!\in\![n]$, 
$q^-_i(t)\!=\!q^\star_i((\alpha+1)\theta-\alpha t)\!=\!0 \Longleftrightarrow 
t\!=\!\theta^{(j)}+\zeta^{(j)}_i$. Furthermore, by construction, 
this common root has multiplicity $1$ for both $q^-_i$ and 
$q^\star_i$. It is then easily checked that 
$h_i(\theta^{(j)})\!=\!\zeta^{(j)}_i$. 

Recalling that the zero scheme of $F$ in $\Ksn$ is exactly 
$\cD_A(F,\cP)\cap\Ksn$ \cite[sec.\ 5.1]{toricint}, we have thus proved 
assertions (1) and (2) of Main Theorem \ref{main:big}. \qed 

\begin{rem} 
The probability of failure in our Las Vegas algorithm above 
is $\boldsymbol{0}$, assuming any probability distribution 
on the coefficients of $F$ yielding probability $1$ avoidance 
of algebraic hypersurfaces in $\K^{\#\mathrm{monomial \ terms}}$. 
\end{rem} 

\subsection{Concluding the Proofs of Main Theorem 1 and 
Corollary \ref{cor:galois}}\mbox{}\\
\label{sub:galois} 
We begin by checking the complexity of our Las Vegas 
algorithm from the preceding section. First note that by Main Theorem  
\ref{main:pert}, each evaluation of $\pert_A$ (for {\bf constant} $u_0,u_1,\ldots,u_n$) takes 
$\cO^*(nR(\oE)^2S(\oE)^{2.376})$ arithmetic steps over $\K$. 
So by observation (a) above (and assertion (0)), we can find $h$ via 
interpolation within time $\cO^*(n\cM(E)R(\oE)^2S(\oE)^{2.376})$. Similarly, by (a), (b), and (c), we can find each $q^-_i$ and $q^\star_i$ within the 
same time bound. So the construction of all these polynomials thus takes a 
total of $\cO^*(n^2\cM(E)R(\oE)^2S(\oE)^{2.376})$ arithmetic steps over $\K$. 

Finding the coefficients of $q^\star_i((\alpha+1)\theta-\alpha t)$ 
takes time $\cO^*(\cM(E)^2)$ via another simple interpolation step. 
So by (d), we can then find $h_1,\ldots,h_n$ still 
within the latter asymptotic time bound. As for space, we only need to 
keep track of $\cO(n\cM(E)^2)$ coefficient values, and this falls 
well within the $\cO(nS(\oE)^2)$ space requirement of Main 
Theorem \ref{main:pert}. 

To conclude, we need only derandomize our algorithm. This 
can be done as follows: replace the generic selection of 
$u_1,\ldots,u_n$ above by $u_i\!=\!\eps^i$ for $i\!\in\![n]$. We 
then obtain that at our genericity condition is violated iff the 
point $(1,\eps,\ldots,\eps^n)\!\in\!\K^{n+1}$   
lies in at least one of $(2n+1)\!$
\begin{tiny}$\begin{pmatrix} \cM(E)\\ 2 \end{pmatrix}$\end{tiny} 
hyperplanes depending on the input $F$. From the box principle, and the 
well-known properties of 
the Van der Monde matrix \cite{binipan}, this can happen to at 
most $n(2n+1)\!$\begin{tiny}$\begin{pmatrix} \cM(E)\\ 2 
\end{pmatrix}$\end{tiny} distinct values of $\eps$. So we can derandomize by 
repeatedly running steps (1)--(3) with new $\eps$ 
at most $n(2n+1)\!$\begin{tiny}$\begin{pmatrix} \cM(E)\\ 2 
\end{pmatrix}$\end{tiny} times, thus finally accounting for 
our aforementioned deterministic time bound. 

Moving on, we must now further refine our algorithm so that our arithmetic 
is over $K$ (or a small algebraic extension thereof) instead of $\K$. 
This can be done as follows: If $\ch \K\!=\!0$, then there are 
enough choices for $\eps$ in $K$ to derandomize our algorithm (since 
$K$ will be infinite). Otherwise, we simply choose $\eps$ in an algebraic 
extension of $K$ of degree\\ 
$\lceil \log_p ((n+1)^2\cM(E)^2)\rceil$, so that we have more than enough 
$\eps$ to choose from. Assertion (3) is now proved. 

To conclude, it immediately follows from   
the development of \cite[sec.\ 6.1]{toricint} that 
the zero scheme of $F$ in $\Kn$ embeds naturally in $\cD_P(F,\cP)$ (and $\cD_\oP(F,\cP)$) if we replace $E$ by $\bO\cup E$. So this introduction 
of extra points into our supports indeed guarantees that $Z'$ includes all 
the affine roots of $F$. \qed 

\begin{rem}
We have thus improved the complexity of finding all the {\bf affine} roots from (roughly) polynomial in $\prod d_i$ to polynomial in $\cM(\bO \cup E)$. 
However, one can improve this even further to polynomial in 
$\sM(E)$ --- the {\bf stable} mixed volume \cite{hsaff,aff} of 
$E$. (In particular, $\sM(E)\!\leq\!\cM(\bO\cup E)\!\leq\!\prod d_i$ 
and the gaps between can be quite large (cf.\ example \ref{ex:spike}).)  
To make this final improvement, it is necessary to use a more refined 
resultant operator --- the {\bf affine} toric resultant, denoted 
$\aff_\oE(\oF)$ \cite{bux}. This is covered at greater length in 
\cite{bux,aff}, and this new operator also allows us to extend Corollaries 
\ref{cor:count} and \ref{cor:double} to $\Kn$ minus an arbitrary union of 
coordinate hyperplanes. 
\end{rem} 

\noindent 
{\bf Proof of Corollary \ref{cor:galois}:} 
It follows immediately from our proof of Main Theorem 1 
that the fields $K[\zeta_i \; | \; (\zeta_1,\ldots,\zeta_n)\!\in\!\Ksn 
\mathrm{ \ is \ a \ root \ of \ } F]$ and $K[\theta \; | \; 
h(\theta)\!=\!0,\prod h_i\!\neq\!0]$ are identical when $u_1,\ldots,u_n$ are 
chosen from $K$. (So the assumption that $\ch \K\!=\!0$ is actually stronger 
than necessary.) Since the latter field is exactly the splitting field of $g$, 
we are done. \qed 

\subsection{The Proof of Main Theorem 2} 
\label{sub:chow} \mbox{}\\ 
We first note that the well known results on the degree 
of $\res_{\bar{E}}(f_1,\ldots,f_{n+1})$ with respect to the 
coefficients of various $f_i$ \cite{combiresult} remain true over 
any algebraically closed field. This follows easily from the 
formulation of the resultant for 
a collection of invertible sheafs on a projective variety \cite{gkz94}. 
In particular, $\choo_A$ should indeed be either be identically zero or a 
homogeneous polynomial (in the $u_a$) of degree $\cM(E)$.  

To prove assertions (1)--(3), we can then simply 
invoke the Vanishing Theorem for Resultants and lemma \ref{lemma:lift}  
(since $\wP$ is compatible with $\conv(A)$). For instance, we obtain  
that $\fii_A(\cZ)$ is positive-dimensional iff $\choo_A$ has infinitely 
many distinct divisors of the form $\sum_{a\in A}\gamma_au_a$. 
So assertion (1) follows immediately. Assertions (2) and 
(3) follow similarly. \qed 

\subsection{The Proof of Corollary \ref{cor:count}}\mbox{}\\ 
Let $\omega$ ($<2.376$) denote the famous matrix multiplication exponent \cite{copper} and set $E_{n+1}\!=\!A$. It then follows immediately 
from \cite{copper} and \cite[The Division Method]{isawres} that for any choice of constant coefficients in $\K$, $\res_\oE(\oF)$ can be evaluated 
within $\cO^*(nR(\oE)S(\oE)^\omega)$ arithmetic operations 
over $\K$, using $\cO(nS(\oE)^2)$ space. 

The first part of Corollary \ref{cor:count} then follows immediately from 
a Van der Monde type argument, as in the proof of Main Theorem 
\ref{main:big}. In particular, via interpolation, it suffices to evaluate $\choo_A(u)$ at exactly $1+n\cM(E)$ distinct points 
of the form $(1,\eps,\ldots,\eps^n)$ to see if 
$\choo_A$ is identically zero. 

To then count the roots of $F$ in $\Ksn$ when 
$\choo_A$ is not identicaly zero, we can begin with a 
variant of the algorithm from Main Theorem 
\ref{main:big} where we evaluate $\choo_{\Delta\cap\Zn}$ instead 
of $\pert_{\Delta\cap\Zn}$. From our previous observations, we 
can thus construct $h$ and $h_1,\ldots,h_n$ within  
time $\cO^*(n^4\cM(E)^3R(\oE)S(\oE)^\omega)$ and 
space $\cO(nS(\oE)^2)$.

We then use the following 
trick: Compute the gcd, $g$, of $h$ and $\prod^n_{i=1} h_i$. 
By remark \ref{rem:mult} of section \ref{sec:chowtor}, we immediately obtain 
that $\deg h - \deg g$ is exactly the number of roots of $F$ in $\Ksn$ 
counting multiplicities. (In fact, the roots of 
$g$ tell us precisely which $\zeta^{(j)}$ lie {\bf out} 
of $\Ksn$.) By the same argument, we can 
also count the number of distinct roots simply by replacing 
$h$ with its square-free part.  By facts (b) and (c) 
of section \ref{sub:solve}, and since the degree of 
$\prod^n_{i=1} h_i$ is at most $n\cM(E)$, these 
computations cause a negligible growth in our asymptotic complexity 
bounds. So we are done. \qed 

\subsection {Facet Searches and the Proof of Main Theorem 3}\mbox{}\\
\label{sub:pfill} 
The first portion of this result follows 
immediately from lemma \ref{lemma:gen} and 
\cite[Corollary 3]{convexapp}. The second portion 
is a consequence of the following algorithm:
\begin{itemize}
\item[{\bf Step 1}]{ Compute the facet normals of 
$P$ and the vertices of all the $\conv(E_i)$. } 
\item[{\bf Step 2}]{ Find a vertex $v$ of some $E_i$ such that for any 
facet normal $w$ of $P$, $v\!\in\!E^w_i \Longrightarrow 
[\#E^w_i\!\geq\!2$ or $\cM(E_1,\ldots,E_{i-1},E_{i+1},\ldots,E_n)\!=\!0]$. 
If no such pair $v$ exists, stop. Otherwise, delete 
any point of $E^w_i$ and go back to step 1. } 
\end{itemize}

By lemma \ref{lemma:gen}, the above algorithm will 
eventually stop with an irreducible fill of $E$. 
As for its complexity, note that the number of facets 
of $P$ is $\cO(m^{2n})$, and we can find the normals to these 
facets within that many arithmetic steps over $\Q$ \cite{gs93}, 
given the convex hulls of the $E_i$. Furthermore, this asymptotic bound dominates the complexity of 
finding the convex hulls of all the $E_i$ \cite{preparata}. 
So the complexity of Step (1) is $\cO(m^{2n})$. Step (2) thus 
amounts to $(n+1)\cO(m^{2n})$ checks for zero mixed volume per 
vertex. So by lemma \ref{lemma:ed}, this takes $\cO(nm^{2n+c})$ arithmetic 
steps over $\Q$. These steps will be executed at most 
$m$ times, so we are done. \qed  

\subsection{Algebraic Homotopies and the Proof of Main Theorem 4}\mbox{}\\
\label{sub:ppert} 
Main Theorem \ref{main:pert} is the cornerstone of our approach 
to solving degenerate systems of equations, so we will precede its 
proof by illustrating one of its underlying constructions: 
explicit algebraic deformation of degenerate zero sets. 

More precisely, following the notation of Main Theorem \ref{main:pert}, 
we will construct a family of curves $C$, fibered over the 
projective line, whose fiber over a particular point  
is a zero-dimensional variety $Z\!\subseteq\!\cZ$ encoding the multiplicities of all the irreducible components of $\cZ$. To do this, we begin with the following lemma, which follows easily from the 
development of \cite[sec.\ 5.1]{toricint} and \cite[sec.\ 11.3]{ifulton}.
\begin{lemma}
\label{lemma:curves} 
Following the notation of definition \ref{dfn:pert} and 
Main Theorem \ref{main:pert}, let $\cZ_0$ be the zero-dimensional 
part of $\cZ$. Also let $\bZ$ be the zero scheme of $F-sF^*$ 
in $\cT(\wP)\times\Pro^1_\K$. Then $\cZ\!=\!\bZ\cap 
(\cT(\wP)\times\{0\})$. Finally, let $C$ be the algebraic curve (possibly reducible) defined by the union of all one-dimensional components 
of $\bZ$ with surjective projection 
onto the second factor of $\cT(\wP)\times\Pro^1_\K$. Then $C$ has the 
following properties:
\begin{enumerate}
\item{ $\bZ\cap(\cT(\wP)\times\{s_0\})\!=\!C\cap(\cT(\wP)\times\{s_0\})$ 
for almost all $s_0\!\in\!\Pro^1_\K$.} 
\item{$Z\!:=\!C\cap(\cT(\wP)\times \{0
\})$ is a subscheme of $\cZ$ consisting of exactly 
$\cM(E)$ points (counting multiplicities). Furthermore, $\cZ_0$ is a 
subscheme of $Z$. } 
\item{Let $W$ be any irreducible component 
of $\cZ$. Then $Z$ has at least one point in $W$ and, for a generic choice of $F^*$, the number of 
points of $Z$ in $W$ (counting multiplicities) is exactly 
the cycle class degree of $W$. \qed } 
\end{enumerate} 
\end{lemma} 
  
We can now begin our most important proof.

\noindent
{\bf Proof of Main Theorem \ref{main:pert}:}  
Similar to the beginning of the proof of Main Theorem \ref{main:chow}, 
the results of \cite{combiresult} (generalized to arbitrary 
algebraically closed $\K$) immediately imply that the degree of 
$\cH$ as a polynomial in $s$ should be $\sum^n_{i=1} 
\cM(E_1,\ldots,E_{i-1},E_{i+1},\ldots,E_n,A)\!\leq\!R(\oE)$.  
Also each coefficient of $\cH(s)$ 
should be a homogeneous polynomial (in the $u_a$) of degree $\cM(E)$. 
These two assertions of course include the opening statement of Main 
Theorem \ref{main:pert} (on the degree and homogeneity of $\pert_A$), but 
they will follow only upon showing that $\cH$ is not identically zero. 

To see this, note that lemma \ref{lemma:once} and the Vanishing 
Theorem for Resultants readily imply that the coefficient of 
the {\bf highest} power of $s$ in $\cH$ is precisely 
$\res_{(E,A)}(F^*,f_{n+1})$. (Simply check the zero set of 
$F-sF^*$ in $\cT(\wP)$ at $s\!=\!\infty$, via 
the homogenization $s'F-sF^*$.) By definition  
\ref{dfn:pert}, and the Vanishing 
Theorem once more, we see that this polynomial in the $u_a$ is not identically 
zero. So $\cH\!\not\equiv\!0$ and we've finished the simplest part 
of our proof. 

Part (1) of Main Theorem \ref{main:pert} follows similarly: One need only 
consider the {\bf unspecialized} resultant polynomial $\res_{(E,A)}(F,f_{n+1})$ 
and observe the terms of degree $0$ in $s$ as we specialize coefficients 
to obtain $F-sF^*$. In particular, $\choo_A(u)$ is precisely 
$\cH(u;0)$. Note then that (2) and (3) also follow almost 
immediately, provided $\choo_A$ is not identically zero. 

To properly handle the cases of (2) and (3) where we are actually 
working with a non-trivial toric perturbation, 
we now invoke lemmata \ref{lemma:lift} and \ref{lemma:curves} to  
establish a precise correspondence between the factors of $\pert_A$ 
and the points of $Z$. 

Letting $\hZ$ be the zero set of $\cH(u;s)$ in $\Pa\times\Pro^1_\K$, 
note that if $k$ is the least exponent of $s$ in $\cH$, 
then $\hZ$ and the zero set of $\frac{\cH}{s^k}$ in 
$\Pa\times\Pro^1_\K$ differ only 
by the presence of the hyperplane $\Pa\times\{0\}$. The second zero 
set does {\bf not} contain this hyperplane, so let's call the second 
zero set $\dZ$. By lemmata \ref{lemma:once} and  
\ref{lemma:lift}, and the Vanishing Theorem for Resultants, we then derive that $\dim[\bZ\cap (\cT(\oP)\times\{s_0\})]\!=\!0$ 
implies the following equivalence:  
$\cH(H_{\fii_A(\zeta)};s_0)\!=\!0 \Longleftrightarrow 
\zeta\!\in\!\bZ\cap (\cT(\wP)\!\times\!\{s_0\})$, where 
$H_p$ is the hyperplane dual to the point $p$.\footnote{So if 
$p\!:=\![p_a \; | \; a\!\in\!A]\!\in\!\Pa$ then $H_p\!:=\!\{[y_a \; | \; 
a\!\in\!A]\!\in\!\Pa \; | \; \sum_{a\in A} p_ay_a\!=\!0\}$.} By 
assertion (1) of lemma \ref{lemma:curves},  
$\dim [\bZ\cap (\cT(\wP)\times\{s_0\})]\!=\!0$ for almost all 
$s_0\!\in\!\Pro^1_\K$. So 
$\fii_A(C)$ is an open subset of $\dZ$, where we define 
$\fii_A(C)\!:=\!\{(y,s_0) \; | \; y\!\in\!H_{\fii_A(\zeta)} \ ; \  
\zeta\!\in\!C\cap(\cT(\wP)\times\{s_0\}) \ ; \ s_0\!\in\!\Pro^1_\K\}$. 
Therefore, since $\fii_A$ is a proper map, $\frac{\cH}{s^k}$ must 
vanish on {\bf all} of $\fii_A(C)$. In particular, via 
remark \ref{rem:mult} of section \ref{sec:chowtor}, 
\[ \pert_A(u)\!=\!\alpha\cdot\!\!\!\prod\limits_{\zeta\in 
C\cap(\cT(\wP)\times\{0\})} 
\left(\sum\limits_{a\in A} \gamma_{\zeta,a}u_a\right) 
\]
where $\alpha\!\in\!\Ks$, $[\gamma_{\zeta,a} \; | \; 
a\!\in\!A]\!:=\!\fii_A(\zeta)$, 
and the product counts intersection multiplicities. 

Continuing our main proof, assertions (2) and (3) follow immediately 
from our last formula and our preceding observations. As for 
the complexity bounds, these follow immediately from our 
earlier fact (a) and the Division Method \cite{isawres} to compute 
$\res_\star(\cdot)$: to evaluate $\pert_A(u)$, we 
simply find the coefficients of $\cH(u;s)$ by 
evaluating $\cH(u;s)$ at $R(\oE)+1$ distinct values of $s$ 
and then interpolating. \qed  

Note that our algebraic proof avoids the use of limiting 
arguments that were present in \cite{gcp}. Thus our result 
holds for any algebraically closed $\K$, instead of just 
$\C$. 

\subsection{Double Perturbations and the Proof of Corollary 
\ref{cor:double}}\mbox{}\\ 
\label{sub:double} 
The first portion of our final corollary follows immediately 
(thanks to Main Theorem \ref{main:pert}) 
by simply replacing $\choo_A$ with $\pert_A$ in the algorithm 
from the proof of Corollary \ref{cor:count}. In particular, 
we obtain that the exact number of roots of $F$ in $\Ksn$ 
(counting multiplicities) is exactly $\deg h^* - \deg g^*$, 
where $h^*$ (resp.\ $g^*$) is the corresponding variant of 
$h$ (resp.\ $g$), using the notation of the proof of Corollary 
\ref{cor:count}. The number of {\bf distinct} roots can of 
course be recovered by using square-free parts (as before), 
thanks to remark \ref{rem:mult} of section \ref{sec:chowtor}. Also, by Main 
Theorem \ref{main:pert}, the complexity of 
this algorithm is just the complexity estimate from 
Corollary \ref{cor:count} multiplied by $R(\oE)$.    

As for the second portion of our corollary, we make a slightly 
more sophisticated variant of the preceding replacement of 
$\choo_A$. 
\begin{dfn}
\label{dfn:double}
Let $F^*$ and $F^{**}$ be $n\!\times\!n$ polynomial 
systems with support contained in $E$ such that (1) $F^*$ and 
$F^{**}$ each have only finitely many roots in 
$\cT(P)$, and (2) $F^*$ and $F^{**}$ share 
no common roots. Following the 
notation of Main Theorem \ref{main:pert}, define a {\bf 
double} toric perturbation of $F$, $\boldsymbol{\pert^{**}_A}$, to be 
the greatest common divisor of $\pert_{A,F^*}$ and $\pert_{A,F^{**}}$. 
\end{dfn}

It is then clear (via Main Theorem \ref{main:pert} once again) that using $\pert^{**}_A$ in place of $\pert_A$ in our preceding algorithm 
will lead to a new estimate, $\deg h^{**} - \deg g^{**}$, for  $\deg(\cZ_0\cap\Ksn)$. Furthermore, by the above definition, it 
is clear that $\deg h^{**}-\deg g^{**}\leq \deg h^* - \deg g^*$. 

As for estimating $\deg \cZ_0$ and $\deg \cZ_\infty$, 
our preceding theory tells us that we can simply respectively 
use $\deg h^{**}$ and $\cM(E)-\deg h^{**}$. \qed
\begin{rem} 
Our algorithm thus requires a generic choice of $F^*$ and $F^{**}$. 
Just as in the construction of $\pert_A$, we can derandomize 
via combinatorial means: We simply use an irreducible fill 
(as in Main Theorem \ref{main:fill}) to construct $F^*$, 
and then simply perturb a single coefficient of $F^*$ 
to construct $F^{**}$. This is the trick used in our earlier 
example in \mbox{section \ref{sub:pert}.}  
\end{rem}  
\begin{rem}
The basic idea behind the double perturbation is that 
the points in $Z\!:=\!\{\gamma(\theta)\}_{h(\theta)=0}$ lying in {\bf 
positive}-dimensional components of $\fii_A(\cZ)$ will move as we vary $F^*$. 
Thus, assuming that $F^{**}$ is such that the new 
$Z$ overlaps the old $Z$ {\bf only} on the isolated roots of $F$,  
we should be able to pick out these isolated roots  
by computing the gcd of $\pert_{A,F^*}$ and $\pert_{A,F^{**}}$.
We hope to address this ``motion of points within a deformation'' 
in future work.  
\end{rem} 
 
\section{Computing Toric Resultants and the Complexity of the Sparse  
Encoding}
\label{sec:complex}
Let us first recall some important facts on the computation of 
toric resultants.

As of 1998, the main method for 
computing $\res_\oE(\oF)$ is to first construct an   
$S(\oE)\times S(\oE)$ {\bf toric resultant  
matrix}, $M_\oE$, whose nonzero entries are certain coefficients 
of $\oF$. This matrix is specifically built so that $\det(M_\oE)$ is,  
for generic choices of the coefficients $c_{i,a}$, a nonzero multiple of $\res_\oE(\oF)$. 
\begin{rem}
So $S(\oE)$ is actually a parameter depending on which algorithm we use for constructing $M_\oE$ --- hence our earlier use 
of an asymptotic bound, instead of an explicit formula, for $S(\oE)$. 
The aforementioned bound is actually a simple estimate on the number of 
lattice points in the interior of the shifted Minkowski sum 
$\delta+\sum^{n+1}_{i=1} E_i$, where $\delta\!\in\!\Qn$ is chosen 
generically. The derivation follows easily from Stirling's estimate 
for the Gamma function, the $n$-dimensional identity 
$\cM(P,\ldots,P)\!=\!n!\V(P)$, and the multilinearity of the 
mixed volume. 
\end{rem}   

Via some clever interpolation tricks \cite{ce,isawres,emipan}, one can 
recover the exact value of $\res_\oE(\oF)$ after interpolating $\det(M_\oE)$ 
through several-many specializations of the coefficients of $\oF$. 
One such fundamental technique, which uses $n+1$ versions of $M_\oE$, 
is known as the {\bf Division Method} \cite{cannyphd,isawres}. 
In general, the matrix $M_\oE$ is highly structured (it is quasi-Toeplitz \cite{emipan}) and, when $\ch\K\!=\!0$, this permits $\res_\oE(\oF)$ to 
be computed much faster than would be expected. 
In practice, the cost of 
building $M_\oE$ (or several versions thereof) can be amortized when one works 
with many $\oF$ with support contained in 
the same $\oE$. Furthermore, when randomization is allowed, 
the results of \cite{ce,isawres,emipan} tell us that this preprocessing is 
actually negligible. 

As for the complexity of computing $\res_\oE(\oF)$ itself, we 
state the following additional facts: 
\begin{enumerate} 
\renewcommand{\theenumi}{\Roman{enumi}}  
\item{\cite[The GCD Method]{isawres} When $\ch \K\!=\!0$, we can 
compute $\res_\oE(\oF)$ (for any choice of constant coefficients in $\K$ 
for $\oF$) within $\cO^*(S(\oE)^{1+\omega})$ arithmetic steps 
and $\cO(S(\oE)^2)$ space.\footnote{ 
The restriction on $\ch \K$ is due to a use of effective 
Hilbert irreducibility, which actually fails in positive 
characteristic \cite{dioplang}.}. However, 
we have the added benefit that we can also compute $\cH(u;s)$ (for 
any {\bf constant} $u\!\in\!\K^{\#A}$) within the {\bf same} 
complexity bound. }
\item{\cite{emipan} If we assume $\ch \K\!=\!0$ {\bf 
and} allow randomization, then we can accelerate the 
Division Method (resp.\ GCD Method) to obtain a Las Vegas time bound of 
$\cO^*(n^2R(\oE)S(\oE)^2)$ (resp.\ $\cO^*(nS(\oE)^3)$). Furthermore, either of 
these improvements requires only $\cO^*(nS(\oE))$ space. }
\item{If $\res_\oE(\oF)\!=\!\det(M_\oE)$, then 
we can reduce the deterministic time bounds of the Division and GCD 
methods to $\cO(R(\oE)^{\omega})$, regardless of $\ch \K$. Furthermore, 
if we also allow randomization and assume $\ch \K\!=\!0$, then we can  
further improve the time bounds of (I) and (II) to $\cO^*(R(\oE)^2)$. 
However, characterizing when $\res_\oE(\oF)$ can be expressed as a 
``small'' determinant is an open problem. (See \cite{wz} for some interesting partial results, including some cases where the Newton 
polytopes are products of scaled standard simplices.) }  
\end{enumerate} 
The last fact is actually a simple corollary of the development 
of \cite{emipan}. In particular, in the situation of (III), 
we can skip an interpolation procedure that would have multiplied our  
time bound by $\cO^*(R(\oE))$. 
 
Let us now state and prove the best current speed-ups for 
all our preceding algorithmic results.  
\begin{cor} 
\label{cor:fast} 
Suppose $\ch \K\!=\!0$ and we allow randomization 
in our algorithms. Then our main algorithmic results 
can be sped up as follows:
\begin{center}
\begin{tabular}{c|c}
 & Sequential (Las Vegas) Time Bound$\ = \cO^*(\cdots)$ \\ \hline 
\vspace{-.4cm} & \\ 
Main Theorem 1 & $n^3\cM(E)R(\oE)^2S(\oE)^2$ or 
                   $n^2\cM(E)S(\oE)^3$ \\
Corollary \ref{cor:count} (First Bound) & $n^2\cM(E)R(\oE)S(\oE)^2$ or 
                   $n\cM(E)S(\oE)^3$ \\ 
Corollary \ref{cor:count} (Second Bound) & $n^3\cM(E)R(\oE)S(\oE)^2$ or 
                   $n^2\cM(E)S(\oE)^3$ \\ 
Main Theorem 4 & $n^2R(\oE)^2S(\oE)^2$ or $nS(\oE)^3$
\end{tabular} 
\end{center} 
Furthermore, the space bound for each of the above 
algorithms is $\cO^*(nS(\oE))$. Finally, if we also have that 
$\res_\oE(\oF)\!=\!\det(M_\oE)$, 
then the four pairs of entries in the right-hand column (from top to 
bottom) can be replaced by the following sequence: $n\cM(E)R(\oE)^3$, 
$\cM(E)R(\oE)^2$, $n\cM(E)R(\oE)^2$, $\cM(E)R(\oE)^3$. \qed 
\end{cor}
\begin{rem} 
As before, the probability of failure in our Las Vegas algorithm above 
is $\boldsymbol{0}$, assuming any probability distribution 
on the coefficients of $F$ yielding probability $1$ avoidance 
of algebraic hypersurfaces in $\K^{\#\mathrm{monomial \ terms}}$. The total 
number of random choices of elements in $K$ (or a small algebraic 
extension thereof) needed is $n+R(\oE)$. (This is just the number choices 
needed to construct $h$ and a variant \cite{emipan} 
of $M_\oE$.)  
\end{rem} 

So in summary, we can solve {\bf any} $n\times n$ system, over an 
algebraically closed field of characteristic zero,  
in Las Vegas time near-quartic in the number of roots of a closely 
related system.\footnote{ We conjecture that this can be done in 
positive characteristic as well. The main current obstruction is 
the use (in current fast algorithms) of algebraic identities for 
recovering elementary symmetric functions from power sums, which 
fail for small positive characteristic. }  Furthermore, 
we can go even faster when we have a sufficiently compressed 
toric resultant matrix. Before proving the above corollary, we will 
briefly explain what we mean by a ``closely related system.'' 

First recall that $\cM(E)$ is precisely the cycle class degree 
of the toric divisor $\cD_P(F,\cP)$ \cite{tfulton,toricint}. Put more simply, 
if we simply perturb the coefficients of $F$, we can expect $F$ 
to have exactly $\cM(E)$ roots in $\cT(P)$. Thus, 
the quantity $\aM_\oE$ defined earlier can be reinterpreted as 
follows: it is the average number of roots of an $n\times n$ 
system of equations with support contained in $(\cE_1,\ldots,\cE_n)$, 
as we let the $\cE_i$ independently range over $\{E_1,\ldots,E_{n+1}\}$, and we 
assume {\bf generically chosen} coefficients. So 
the quantity $S(\oE)$ can also be interpreted as a weighted average of a set 
of cycle class degrees. Similarly, note that the 
generic number of roots of the $(n+1)\times (n+1)$ 
system $(F-sF^*,s-s_0)$ is exactly $\cM(E_1\times \{0,1\},\ldots,
E_n\times \{0,1\},\{0,\hat{e}_{n+1}\})$. So by the 
multilinearity of the mixed volume, the last mixed volume is exactly 
$R(\oE)$. 

Let us now prove our above corollary.

\noindent
{\bf Proof of Corollary \ref{cor:fast}:} 
The key bounds to begin with are those of Corollary \ref{cor:count}. In 
particular, the first bound of Corollary \ref{cor:count} is the complexity of determining whether 
$\choo_A(u)$ vanishes identically. Since this can be accomplished 
by evaluating $\res_\oE(\oF)$ at $\cM(E)+1$ random points, 
facts (I)--(III) above immediately imply our asserted bounds.

As for the second bound of Corollary \ref{cor:count}, this is the complexity 
of running a variant of the algorithm of Main Theorem 1, where 
$\pert_A$ is replaced by $\choo_A$. Since $\choo_A$ is just a 
specialized resulant, and since this algorithm boils down 
to evaluating $\choo_A$ at $\cO(n\cM(E))$ distinct points, facts (I)--(III) 
immediately imply these bounds as well.

From the proof of Main Theorem 4, we know that 
the bound from Main Theorem 4 is simply the complexity of 
evaluating $\res_\oE(\oF)$ at $<\!R(\oE)\!+\!1$ different specializations 
of $s$. (Remember that $s$ occurs only in the coefficients of  
$f_1,\ldots,f_n$, and all other parameters are assumed to be constants.) 
So this bound follows easily from facts (I)--(III) as well. 

To conclude, the bound from Main Theorem 1 is simply the 
complexity of evaluating $\pert_A$ at $\cO(n\cM(E))$ 
distinct points. From the bound of Main Theorem 4, 
we are done. \qed  

Are the above complexity bounds the best one can expect for solving 
polynomial systems  
specified in the sparse encoding?\footnote{That is, when we specify polynomial 
systems as a list of exponents and coefficients.} Neglecting the precise 
values of the exponents (which we've seen range somewhere 
between $4$ and $7.376$, if not better), the answer is ``yes.'' 
This is due to the fact that a generic $F$ will have exactly $\cM(E)$ 
distinct roots in $\K$, {\bf regardless} of the number of terms present. Thus, 
it is really $\cM(E)$, {\bf not} the number of terms, which 
governs the complexity of {\bf global} polynomial system solving 
over an algebraically closed field. So the quantities in the 
``base'' of our bounds can not be any smaller (asymptotically) than 
$\cM(E)$. As for the exponent, we so far have only the obvious 
worst case lower bound of $1$. 

However, the question of whether the number of terms more strongly 
governs the complexity of solving over a {\bf non}-algebraically 
closed field, or solving for a {\bf single} root, is quite 
open for $n\!\geq\!2$. For example, while Khovanskii has shown that 
the number of {\bf real} roots of a sparse system of equations is singly 
exponential in 
the number of terms \cite{sparse}, the complexity of real solving  
is {\bf not} yet known to fall within such a bound if $n\!\geq\!2$. 
Similarly, 
while a recent algorithm of Ye \cite{ye} for $\eps$-approximating a 
{\bf single} $d^{\thth}$ root of $\alpha\!\in\!\R$ has complexity 
$\cO((\log d)\log\log\frac{|\alpha|}{\eps})$, the complexity of 
finding a {\bf single} root of $F$ in $\Kn$ is quite open. 
We hope to address these finer points of sparse algebraic complexity 
in future work.  
 
We now close with a brief example of how $\cM(E)$ can 
be smaller than $d_\Pi$ (the product of the total degrees of $f_1,\ldots,f_n$) by an exponential factor.
\begin{ex} {\bf (Well Directed Spikes)} 
\label{ex:spike}
Consider the system of equations $F$ defined by 
\[ a_{1,1}+a_{1,2}x_1+\cdots+a_{1,n}x_{n-1}+c_{1,1}(x_1\cdots x_n) 
+ \cdots + c_{1,d}(x_1\cdots x_n)^d =0 \]
\[ \vdots \] 
 \[ a_{n,1}+a_{n,2}x_1+\cdots+a_{n,n}x_{n-1}+c_{n,1}(x_1\cdots x_n) 
+ \cdots + c_{n,d}(x_1\cdots x_n)^d =0. \] 
In this case, the Newton polytopes are all equal to a single  
``spike,'' and this spike is equivalent (via an integer linear 
map with determinant $1$) to a standard $n$-simplex 
scaled by $d$ in the $x_1$-direction. So it is easy 
to check that $\cM(E)\!=\!d$. However, the product of 
the total degrees of $F$ is clearly $n^nd^n$. (It 
is also not hard to see that the best {\bf multigraded} B\'ezout bound 
\cite{wampler} is $n!d^n$.) Generating infinite families of such examples is 
easy, simply by picking Newton polytopes which are $n$-dimensional, 
but ``long'' in a suitable fixed direction.  
\end{ex} 

\begin{rem}
The construction of toric resultant matrices is an area of 
active research and it can be reasonably expected that our 
earlier asymptotic estimate on $S(\oE)$ will be significantly 
improved in the near future. In particular, a significant 
first step would be to find an algorithm which always constructs a   
toric resultant matrix of size $\cO(R(\oE))$. Looking even 
further ahead, there is also hope for general algorithms which 
construct even smaller matrices, via the use of entries which 
are nonlinear polynomials in the coefficients of $\oF$.  
\end{rem}   

\section{Acknowledgements}
The author would like to thank an anonymous referee 
for extensive comments on clarifying 
the exposition and development of this paper. The author 
also expresses his deep gratitude to Gregorio Malajoich-Mu\~noz 
for his assistance in computing the matrix from section  
\ref{sub:pert}. Special thanks also go to Pino Italiano and Nini Wong 
for their help in obtaining a copy of \cite{pspace}. 

This paper is dedicated to my son, Victor Lorenzo.  

\bibliographystyle{amsalpha}

\begin{thebibliography}{A}

\bibitem[Ber75]{bernie} Bernshtein, D. N., {\it ``The Number of 
Roots of a System of Equations,"} Functional Analysis and its Applications (translated from Russian), Vol. 9, No. 2, (1975), 
pp.\ 183--185.

\bibitem[BP94]{binipan} Bini, Dario and Pan, Victor Y. {\it Polynomial and
Matrix Computations, Volume 1: Fundamental Algorithms,} Progress in
Theoretical Computer Science, Birkh\"auser, 1994.

\bibitem[BZ88]{buza} Burago, Yu. D. and Zalgaller, V. A., {\it
Geometric Inequalities,} Grundlehren der mathematischen Wissenschaften 285,
Springer-Verlag (1988).

\bibitem[Can87]{cannyphd} Canny, John F., {\it ``The Complexity of 
Robot Motion Planning Problems,''} ACM Doctoral Dissertation Award 
Series, ACM Press (1987). 

\bibitem[Can88]{pspace} \underline{\hspace{\jfc}}, {\it ``Some Algebraic 
and Geometric Computations in PSPACE,''} Proc.\ 20$^\thth$ ACM 
Symp.\ Theory of Computing, Chicago (1988). 

\bibitem[Can90]{gcp} \underline{\hspace{\jfc}}, {\it ``Generalised 
Characteristic Polynomials,''} J.\ Symbolic Computation 9 (1990), 
no.\ 3, pp.\ 241--250. 

\bibitem[CE95]{ce} Canny, John F. and Emiris, Ioannis
Z.,  {\em ``A Subdivision-Based Algorithm for
the Sparse Mixed Resultant,''} preprint, INRIA (1995).

\bibitem[CKL89]{ckl89} Canny, J. F., Kaltofen, Eric, and 
Lakshman, Y., {\it 
``Solving Systems of Non-Linear Polynomial Equations Faster,''} 
Proc.\ ACM Intern.\ Symp.\ on Symbolic and Algebraic Computation, 
pp.\ 121--128, 1989. 

\bibitem[CG84]{chigo} Chistov, A. L., and Grigoriev, Dima Yu, {\it 
``Complexity of Quantifier Elimination in the Theory of Algebraically 
Closed Fields,''} Lect.\ Notes Comp.\ Sci.\ 176, Springer-Verlag (1984). 

\bibitem[CW90]{copper} Coppersmith, Don and Winograd, Shmuel, 
{\it ``Matrix Multiplication via Arithmetic Progressions,''} J.\ Symbolic Computation, 9 (1990), no.\ 3, pp.\ 251--280. 

\bibitem[Dan78]{dannie} Danilov, V. I., {\it ``The Geometry of Toric
Varieties,"} Russian Mathematical Surveys, 33 (2), pp.\ 97--154, 1978.

\bibitem[DK95]{black} D\'\i{}az, Angel, and Kaltofen, 
Erich, {\it ``On Computing Greatest Common Divisors 
with Polynomials Given by Black Boxes for Their 
Evaulations,"} Proceedings of ISSAC '95, Montreal 
Canada, pp.\ 232--239, ACM Press (1995).

\bibitem[DGH98]{mvcomplex} Dyer, M., Gritzmann, P., and 
Hufnagel, A., {\it ``On the Complexity of Computing Mixed 
Volumes,"} SIAM J.\ Comput. {\bf 27} (1998), no.\ 2, 
356--400.

\bibitem[DS95]{introchow} Dalbec, John, and Sturmfels, Bernd, 
{\it ``Introduction to Chow Forms,"} Invariant Methods in Discrete 
and Computational Geometry (Cura\c{c}ao, 1994), pp.\ 37--58, Kluwer Academic 
Publishers, Dordrecht, 1995.  

\bibitem[EC95]{isawres} Emiris, Ioannis Z. and Canny, John F., {\it
``Efficient Incremental Algorithms for the Sparse Resultant and the Mixed
Volume,''} Journal of Symbolic Computation, vol.\ 20 (1995), pp.\ 117--149.

\bibitem[EP97]{emipan} Emiris, Ioannis Z. and Pan, Victor 
Y., {\it ``The Structure of Sparse Resultant Matrices,''} 
Proceedings of the International Symposium on Symbolic 
and Algebraic Computation (ISSAC) 1997, ACM Press.

\bibitem[Ewa96]{ewald} Ewald, G\"unter, {\it Combinatorial 
Convexity and Algebraic Geometry,} Graduate Texts in Mathematics 
168, Springer-Verlag, New York, 1996. 

\bibitem[Ful84]{ifulton} Fulton, William, {\it Intersection Theory,} Springer-Verlag, 1984.

\bibitem[Ful93]{tfulton} \underline{\hspace{\fw}}, {\it Introduction to Toric
Varieties}, Annals of Mathematics Studies, no.\ 131, Princeton University
Press, Princeton, New Jersey, 1993.

\bibitem[GKZ94]{gkz94} Gel'fand, I. M., Kapranov, M. M., and 
Zelevinsky, A. V., 
{\it Discriminants, Resultants and Multidimensional Determinants,} 
Birkh\"auser, Boston, 1994. 

\bibitem[GHMP95]{joos} Giusti, M., Heintz, J., Morais, 
J.\ E., Pardo, L.\ M., {\it ``When Polynomial Equation 
Systems can be `Solved' Fast?,"} Applied Algebra, Algebraic 
Algorithms and Error-Correcting Codes (Paris, 1995), 
205--231, Lecture Notes in Comput.\ Sci.\ 948, Springer, 
Berlin, 1995. 

\bibitem[GV91]{lalo} Gonz\'alez-Vega, Laureano,  {\it 
``A Subresultant Theory for Multivariate Polynomials,''} 
Proceedings of the 1991 International Symposium on Symbolic and 
Algebraic Computation, pp.\ 79--85, Stephen M. Watt (ed.), ACM Press. 

\bibitem[GK94]{gk94} Gritzmann, Peter and Klee, Victor, 
{\it ``On the Complexity of Some Basic Problems in Computational 
Convexity II: Volume and Mixed Volumes,''} Polytopes: Abstract, 
Convex, and Computational (Scarborough, ON, 1993), pp.\ 
373--466, NATO Adv.\ Sci.\ Inst.\ Ser.\ C Math.\ Phys.\ Sci., 
440, Kluwer Acad.\ Publ., Dordrecht, 1994. 

\bibitem[GS93]{gs93} Gritzmann, Peter and Sturmfels, Bernd, 
{\it ``Minkowski Addition of Polytopes: Computational Complexity and Applications to Grobner Bases,''} SIAM J.\ Discrete Math.\ 6 (1993), no.\ 
2, 246--269.  

\bibitem[HS97]{hsaff} Huber, Birk and Sturmfels, Bernd, {\it ``Bernshtein's 
Theorem in Affine Space,"} Discrete and Computational Geometry, {\bf 17} 
(1997), no.\ 2, 137--141. 

\bibitem[KSZ92]{ksz92} Kapranov, M. M., Sturmfels, B., and Zelevinsky, 
A. V., {\it ``Chow Polytopes and General Resultants,''} 
Duke Mathematical Journal, Vol.\ 67, No.\ 1, July, 1992, pp.\ 189--218. 

\bibitem[KKMS73]{kkms} Kempf, G., Knudsen, F., Mumford, D., Saint-Donat,
B., {\it Toroidal Embeddings I}, Lecture Notes in Mathematics 339,
Springer-Verlag, 1973.

\bibitem[Kho77]{khocompat} Khovanskii, Askold G., {\it ``Newton 
Polyhedra and Toroidal Varieties,"} Functional Anal.\ Appl., 11 (1977), 
pp.\ 289--296.

\bibitem[Kho91]{sparse} \underline{\hspace{\agk}}, {\it Fewnomials,}
AMS Press, Providence, Rhode Island, 1991. 

\bibitem[Lan83]{dioplang} Lang, Serge, {\it Fundamentals of 
Diophantine Geometry,} Springer, New York, 1983. 

\bibitem[Mal94]{gregthesis} Malajovich-Mu\~noz, Gregorio, {\it ``The 
Complexity of Newton Iteration and Path-Following Algorithms,''} 
Ph.D.\ thesis, University of California, available from University 
Microfilms International, Michigan (1994). 

\bibitem[MM82]{mm} Mayr, E. and Meyer, A., {\it ``The 
Complexity of the Word Problem for Commutative Semigroups 
and Polynomial Ideals,''} Adv.\ Math.\ {\bf 46}, 305--329, 1982. 

\bibitem[MP98]{mp98} Mourrain, Bernard and Pan, Victor Y. {\it 
``Asymptotic Acceleration of Solving Multivariate Polynomial 
Systems of Equations,''} Proc.\ ACM STOC 1998.   

\bibitem[PS93]{chowprod} Pedersen, P. and Sturmfels, B., {\it ``Product
Formulas for Sparse Resultants and Chow Forms,"} Mathematische Zeitschrift,
214: 377--396, 1993.

\bibitem[PS85]{preparata} Preparata, Franco P. and
Shamos, Michael Ian, {\it Computational Geometry: An Introduction,}
Texts and Monographs in Computer Science, Springer-Verlag,
New York-Berlin, 1985.

\bibitem[Ren87]{renegar} Renegar, Jim, {\it ``On the Worst Case 
Arithmetic Complexity of Approximating Zeros of Systems of Polynomials,''} 
Technical Report, School of Operations Research and Industrial Engineering, 
Cornell University.

\bibitem[Roj94]{convexapp} Rojas, J. Maurice, {\it ``A Convex Geometric
Approach to Counting the Roots of a Polynomial System,"} Theoretical
Computer Science (1994), vol.\ 133 (1), pp.\ 105--140.

\bibitem[Roj97a]{rio} \underline{\hspace{\jmr}}, {\it ``Toric 
Laminations, Sparse Generalized Characteristic Polynomials, and a Refinement 
of Hilbert's Tenth Problem,''} Foundations of Computational Mathematics,
selected papers of a conference, held at IMPA in Rio de Janeiro, January
1997, Springer-Verlag (1997).   

\bibitem[Roj97b]{bux} \underline{\hspace{\jmr}}, {\it 
``Affine Elimination Theory,''}  
extended abstract, Proceedings of a Conference in Honor 
of the 60$\thth$ birthday of David A.\ Buchsbaum, Northeastern 
University, October, 1997. 

\bibitem[Roj98a]{venice} \underline{\hspace{\jmr}}, {\it ``Intrinsic  
Near Quadratic Complexity Bounds for Real Multivariate Root Counting,''} 
Proceedings of the Sixth Annual European Symposium on Algorithms, 
Lecture Notes in Computer Science 1461, Springer-Verlag (1998). 

\bibitem[Roj98b]{resvan} \underline{\hspace{\jmr}}, {\it 
``The Geometry of Elimination I: Degree Formulae and the 
Vanishing of Resultants,'' } preprint, City University of Hong 
Kong (1998). 

\bibitem[Roj98c]{aff} \underline{\hspace{\jmr}}, {\it 
``The Geometry of Elimination II: Affine Elimination 
Theory and Better Nullstellensatze,'' } preprint, City University 
of Hong Kong (1998). 

\bibitem[Roj99]{toricint} \underline{\hspace{\jmr}}, {\it ``Toric Intersection 
Theory for Affine Root Counting,''} Journal of Pure and 
Applied Algebra, June 1999. 

\bibitem[RW96]{rojaswang} Rojas, J. M., and Wang, Xiaoshen, {\it 
``Counting Affine Roots of Polynomial Systems Via Pointed Newton 
Polytopes,"} Journal of Complexity, vol.\ 12, June (1996), pp.\ 
116--133.

\bibitem[Sch94]{schneider} Schneider, Rolf, {\it Convex Bodies: The
Brunn-Minkowski Theory,} Encyclopedia of Mathematics and its Applications,
v.\ 44, Cambridge University Press, 1994. 

\bibitem[Stu93]{sparseelim} Sturmfels, Bernd, {\it ``Sparse Elimination 
Theory,''} In D. Eisenbud and L. Robbiano, editors, Proc. 
Computat.\ Algebraic Geom.\ and Commut.\ Algebra 1991, pages 
377--396, Cortona, Italy, 1993, Cambridge Univ. Press. 

\bibitem[Stu94]{combiresult} \underline{\hspace{\bernd}}, {\it ``On the Newton 
Polytope of the Resultant,''} Journal of Algebraic Combinatorics, 3: 207--236, 
1994. 

\bibitem[Stu98]{intro} \underline{\hspace{\bernd}}, {\it ``Introduction to
Resultants,''} Applications of Computational Algebraic Geometry 
(San Diego, CA, 1997), 25--39, Proc.\ Sympos.\ 
Appl.\ Math., 53, Amer.\ Math.\ Soc., Providence, RI, 1998.

\bibitem[Van50]{vdv} van der Waerden, B. L., {\it Modern Algebra,} 
2$^\nd$ edition, F.\ Ungar, New York, 1950.

\bibitem[Wam92]{wampler} Wampler, Charles W., {\it ``Bezout Number
Calculations for Multi-Homogeneous Polynomial Systems,''} Applied
Mathematics and Computation 51, pp.\ 143--157, 1992.

\bibitem[WZ94]{wz} Weyman, Jerzy and Zelevinsky, Andrei, {\it ``Determinantal Formulas for Multigraded Resultants,''} 
J.\ Algebraic Geom.\ 3 (1994), no.\ 4, 569--597.

\bibitem[Ye94]{ye} Ye, Yinyun, {\it ``Combining Binary Search 
and Newton's Method to Compute Real Roots for a Class of 
Real Functions,''} J. Complexity 10 (1994), no.\ 3, 271--280.

\end{thebibliography}

\end{document}